\documentclass[onefignum,onetabnum]{siamonline220329}
\usepackage{import}

\usepackage{mathrsfs}
\usepackage{mathtools}
\usepackage{amssymb}

\usepackage{comment}
\usepackage{soul}

\crefname{Assumption}{Assumption}{Assumptions}

\usepackage{cleveref}

\newcommand{\br}{{\normalfont{\boldsymbol{r}}}}
\newcommand{\bxi}{{\normalfont{\boldsymbol{\xi}}}}

\newcommand{\dd}{{\normalfont{\text{d}}}}

\usepackage{xcolor}
\definecolor{ao}{rgb}{0.0, 0.5, 0.0}

\DeclareSymbolFont{extraup}{U}{zavm}{m}{n}
\DeclareMathSymbol{\varheart}{\mathalpha}{extraup}{86}
\DeclareMathSymbol{\vardiamond}{\mathalpha}{extraup}{87}

\usepackage{tikz}

\usepackage{pifont}

\usepackage{color}
\newcommand{\half}{\frac{1}{2}}

\newcommand{\norm}[1]{\left \lVert #1 \right \rVert}
\newcommand{\snorm}[1]{\left \lvert #1 \right \rvert}

\newcommand{\IR}{\mathbb{R}}
\newcommand{\IC}{\mathbb{C}}

\newcommand{\IN}{\mathbb{N}}

\newcommand{\OT}{\mathsf{T}}

\newcommand{{\D}}{\normalfont{\text{D}}}
\newcommand{{\G}}{\normalfont{\text{G}}}
\newcommand{{\U}}{\normalfont{\text{U}}}
\newcommand{{\I}}{\normalfont{\text{I}}}
\newcommand{{\N}}{\mathsf{N}}
\newcommand{{\per}}{\normalfont{\text{per}}}
\newcommand{{\y}}{{\boldsymbol{{y}}}}

\newcommand{{\bc}}{{\bf c}}
\newcommand{{\bx}}{{\bf x}}
\newcommand{{\by}}{{\bf y}}
\newcommand{{\bz}}{{\bf z}}
\newcommand{{\bd}}{{\bf d}}

\newtheorem{problem}[theorem]{Problem}

\newcommand{\dual}[2]{\left \langle #1,#2 \right \rangle}

\newcommand{\uinc}{u^{\normalfont\text{inc}}}

\usepackage{hyperref}
\usepackage{cleveref}
\usepackage{autonum}

\usepackage{graphicx}
\usepackage{subcaption}
\usepackage{caption}

\usepackage{lipsum}
\usepackage{amsfonts}
\usepackage{graphicx}
\usepackage{epstopdf}
\usepackage{algorithmic}
\ifpdf
  \DeclareGraphicsExtensions{.eps,.pdf,.png,.jpg}
\else
  \DeclareGraphicsExtensions{.eps}
\fi

\usepackage{enumitem}
\setlist[enumerate]{leftmargin=.5in}
\setlist[itemize]{leftmargin=.5in}

\newsiamremark{remark}{Remark}
\newsiamremark{assumption}{Assumption}
\newsiamremark{hypothesis}{Hypothesis}
\crefname{hypothesis}{Hypothesis}{Hypotheses}
\crefname{assumption}{assumption}{Assumptions}
\newsiamthm{claim}{Claim}

\headers{Domain UQ for the Lippmann-Schwinger
VIE}{F. Henr\'iquez and I. Labarca-Figueroa}

\title{Domain Uncertainty Quantification for the Lippmann-Schwinger
Volume Integral Equation\thanks{Version of \today.
\funding{The work of IL-F was funded by the Austrian Science Fund (FWF) P 35673-N.}}}

\author{
Fernando Henr\'iquez\thanks{Chair of Computational Mathematics and Simulation Science (MCSS), \'Ecole Polytechnique F\'ed\'erale de Lausanne, Lausanne, Switzerland
(\email{fernando.henriquez@epfl.ch}).}
\and 
Ignacio Labarca-Figueroa\thanks{Institute for Theoretical Physics, University of Innsbruck, Innsbruck, Austria 
(\email{ignacio.labarca-figueroa@uibk.ac.at}).}
}

\usepackage{amsopn}

\begin{document}

\maketitle

\begin{abstract}
In this work, we consider the propagation of acoustic waves in unbounded domains characterized by a constant wavenumber, except possibly in a bounded region. The geometry of this inhomogeneity is assumed to be uncertain, and we are particularly interested in studying the propagation of this behavior throughout the physical model considered. A key step in our analysis consists of recasting the physical model—originally set in an unbounded domain—into a computationally manageable formulation based on Volume Integral Equations (VIEs), particularly the Lippmann-Schwinger equation. We show that both the leading operator in this volume integral formulation and its solution depend holomorphically on shape variations of the support of the aforementioned inhomogeneity. This property, known as shape holomorphy, is crucial in the analysis and implementation of various methods used in computational Uncertainty Quantification (UQ). We explore the implications of this result in forward and inverse UQ and provide numerical experiments illustrating and confirming the theoretical predictions.
\end{abstract}

\begin{keywords}
Volume Integral Equations, Shape Uncertainty Quantification, Higher-Order Quasi-Monte Carlo, Bayesian Inverse Problems, Shape Holomorphy.
\end{keywords}

\begin{MSCcodes}
32D05, 35A20, 45P05, 65R20.
\end{MSCcodes}

%
%
%

\section{Introduction}\label{sec:intro}	
We consider the acoustic scattering at a penetrable, bounded, inhomogeneous domain of uncertain shape. We are interested in studying the propagation of uncertainty due to domain variations in this model and their effect on the so-called Quantities of Interest (QoI). The quantification of the effect of shape deformations has become an increasingly important subject due to its significant impact in various fields of science and engineering. Unlike previous work in this subject, we adopt a new viewpoint for the mathematical model describing the underlying wave propagation phenomena. This approach involves using VIEs to simultaneously and consistently address the unbounded nature of the problem and the presence of bounded, possibly not piecewise constant, inhomogeneous inclusions. Well-posedness of VIEs under different assumptions on material properties is studied in \cite{Cos88}. Fast solvers for the Lippmann-Schwinger integral equation based on collocation schemes are presented in \cite{ambikasaran2016fast,bruno2004efficient,vainikko2000fast, anderson2024fast}. In this work, we study Galerkin formulations of the VIEs due to their robustness and provable stability \cite{kirsch2009operator},\cite[Chapter~1]{labarca2024coupled}. Both properties will show to be important in the context of parametric domain deformations
studied in this work.

A variety of methods have been proposed to computationally handle domain uncertainty quantification.
One approach is to perform a perturbation analysis under the assumption of small domain deformations with respect to a nominal shape, as in \cite{harbrecht2008sparse,harbrecht2013first,harbrecht2015second}.
This method relies on shape calculus techniques and suitable linearization with respect to the nominal shape. Unfortunately, this approach becomes less accurate as the amplitude of domain deformations increases and we leave the small deformations regime.

The \emph{domain mapping approach} has been introduced in \cite{xiu2006numerical, tartakovsky2006stochastic} to properly account for large domain deformations. Therein, as the name suggests, parameter-dependent domain transformations defined in a reference domain are used to represent each possible domain deformation. This approach has been further explored in \cite{CNT2016} and \cite{HPS16} with a focus on elliptic PDEs in random domains, in \cite{AJS19} for electromagnetic wave phenomena, in \cite{harbrecht2024quantifying} for linear elasticity, and in \cite{DHJM2022} for acoustic wave propagation in unbounded domains. In particular, the latter study used boundary integral operators and a suitable boundary integral formulation for sound-soft acoustic scattering.

A commonality of the previously described works is that in the definition of the 
domain deformations parametrically-defined, high-dimensional shape variations
are considered. This setting naturally arises when the parametric 
shape representations are obtained as the Karhunen-Lo\`{e}ve expansion of a random field.
This, in turn, defines a parameter-to-solution map or parameter-to-QoI map with a high-dimensional,
possibly countable infinite, parametric input.

A variety of techniques and methods have been proposed to computationally 
handle the efficient approximation of parametric maps with high-dimensional inputs.
Examples of these include sparse grid interpolation and quadrature \cite{ZS17,SS13,NTW2008,HHPS2018},
higher-order Quasi-Monte Carlo integration (HoQMC) \cite{DKL14,DLC16,DGLS19,DGLS17},
construction of neural network-based surrogates \cite{SZ19,HSZ20,OSZ22,HOS22,ABD22,weder2024galerkin},
and the reduced basis method \cite{chen2016model}.
As pointed out in \cite{CCS15}, a key property to break the so-called \emph{curse of dimensionality}
in the parameter space is the holomorphic dependence of the parameter-to-solution map.
This property has been established for a variety of problems, including for example 
subsurface flows \cite{CNT2016,HPS16,HS2022}, time-harmonic electromagnetic wave
scattering \cite{AJS19}, stationary Stokes
and Navier-Stokes equation \cite{CSZ18}, Helmholtz equation
\cite{HSSS15,SW23,GKS21}, linear elasticity \cite{harbrecht2024quantifying}, and for boundary 
integral operators \cite{HS21,henriquez2021shape,PHJ23,DH23,DL23,riva2022shape}.

\subsection*{Contributions}
In this work, we explore the application of various computational techniques to the forward and inverse shape UQ of the acoustic scattering problem involving a compactly supported inhomogeneous inclusion. Our approach is based on the formulation of the exterior scattering problem as an equivalent VIE. Unlike the approach followed in \cite{HSSS15}, we do not need to include an artificial boundary to truncate the computational domain. Furthermore, we do not need to construct domain shape deformations that become the identity as we approach said artificial boundary. On the other hand, VIEs allow us to include material inhomogeneities that are not piecewise constant, which is a limitation of the approach proposed in \cite{DHJM2022}.

As a stepping stone in the analysis of our problem, we prove that the Volume Integral Operators (VIOs) appearing in our formulation depend holomorphically on shape deformations. Among the consequences of this property are that the so-called domain-to-solution map is holomorphic, as is its discrete counterpart obtained by Galerkin solutions of VIEs, and so is any quantity of interest (QoI) defined upon them. This result enables us to use technical results for the approximation of maps with high-dimensional parametric inputs.


\subsection*{Outline}
This work is organized as follows. \Cref{sec:preliminaries} is devoted
to introducing relevant concepts and results to be used throughout this work.
In particular, in \Cref{sec:holomorphy_banach} we review the notion of holomorphy in Banach 
spaces, and in \Cref{sec:rt_theorem} we introduce the Riesz-Thorin 
interpolation theorem, which is of key importance to establish boundedness of
certain boundary integral operators.
In \Cref{sec:vol_int_eq} we consider
the Helmholtz equation with compactly supported variable coefficients, and explain
its reformulation in terms of VIEs. 
Following the works \cite{LH23,labarca2024coupled}, we recall well-posedness of the continuous formulation
and its discrete counterpart.
In \Cref{sec:shape_holomorphy} we present a thorough
shape holomorphy analysis of involved VIOs and of the solution to the 
corresponding VIEs.
In \Cref{sec:applications} we consider two applications related to forward and inverse UQ
and discuss the significance of the previously established shape holomorphy result in 
this context. 
In \Cref{sec:numerical_experiments} we present a set of numerical 
experiments involving the approximation of a given QoI in forward and inverse UQ.
We conclude this work in \Cref{sec:concluding_remarks} providing final observations
and sketching directions of future research.

\section{Preliminaries}
\label{sec:preliminaries}

\subsection{Notation}
\label{sec:notation}
Let $\mathcal{O} \subset \mathbb{R}^d$, $d \in \{2,3\}$, be a bounded Lipschitz domain
with boundary $\partial \mathcal{O}$. For $s\geq 0$, by $H^s(\mathcal{O})$ we denote the standard Sobolev
spaces in $\mathcal{O}$, with $L^2(\mathcal{O}) \equiv H^0(\mathcal{O})$ and by $\widetilde{H}^{-s}(\mathcal{O})$ its dual with respect to the $L^2(\mathcal{O})$-based duality paring, which we denote by $\dual{\cdot}{\cdot}_{\mathcal{O}}$.
For (complex) Banach spaces $X$ and $Y$,
we denote by $\mathscr{L}(X,Y)$ the space of bounded linear operators 
from $X$ into $Y$ and by $\mathscr{L}_{\text{iso}}(X,Y)$ the (open)
subset of isomorphisms, i.e.~bounded linear operators with a bounded inverse. 
Recall that $\mathscr{L}(X,Y)$
is a (complex) Banach space equipped with the standard operator norm \cite[Theorem III.2]{RS80_Vol1}.

Finally, we set $\mathbb{U}=[-1,1]^{\IN}$ and equip it with the product topology.
According to Tychonoff's theorem, this renders $\mathbb{U}$ compact with this topology.
\subsection{Holomorphy in Banach Spaces}
\label{sec:holomorphy_banach}
The main theoretical result introduced in this paper corresponds to the analytic dependence
of VIEs upon domain deformations. To this end, we introduce the precise 
notion of holomorphy in Banach spaces to be used throughout this work. 

\begin{definition}[{\cite[Definition 13.1]{Muj86}}]\label{def:complex_diff}
Let $U$ be an open, nonempty subset of $E$. 
A map {$\mathcal{F}:U \rightarrow F$} is said to be
\emph{Fr\'echet differentiable}
if for each $r\in U$ there exists a map
$(\frac{d}{d r} \mathcal{F})(r,\cdot)\in \mathscr{L}(E,F)$ such that
\begin{align}
	\norm{
	{{\mathcal{F}(r+\xi)}}
	-
	{
	\mathcal{F}(r)
	}
	- 
	\left(\frac{d}{dr} \mathcal{F}\right)(r,\xi)}_F 
	=
	o\left(\norm{\xi}_E\right).
\end{align}
We say that {$\left(\frac{d}{dr} \mathcal{F}\right)(r,\xi)$}
is the \emph{Fr\'echet derivative} of the map
{$\mathcal{F}: U \rightarrow F$} at $r\in U$ in the
direction $\xi\in E$.
\end{definition}

Recursively, one may define higher-order Fr\'echet derivatives.

\begin{definition}
Let $U$ be an open, nonempty subset of $E$ and let $m\in \IN$. 
We say that the map {$\mathcal{F}:U \rightarrow F$} is
$m$-times Fr\'echet differentiable if it is $(m-1)$-times
Fr\'echet differentiable and the map
{
\begin{align}
	U \ni
	r
	\mapsto
	\left(\frac{d^{m-1}}{dr^{m-1}} \mathcal{F}\right)
	(r,\underbrace{\cdot,\dots,\cdot}_{m-1 \text{ times}})
	\in
	\mathscr{L}\left(E^{(m-1)},F\right)
\end{align}
}%
is Fr\'echet differentiable as well. 
We say that ${\mathcal{F}:U \rightarrow F}$ is infinitely complex
Fr\'echet differentiable if it is $m$-times Fr\'echet differentiable for all
$m\in \mathbb{N}$.
\end{definition}

In the following, we adopt the notation 
\begin{align}
	\left(\frac{d^{m}}{dr^{m}} \mathcal{F}\right)(r,\xi)
	=
	\left(\frac{d^{m}}{dr^{m}} \mathcal{F}\right)
	(r,\underbrace{\xi,\cdots,\xi}_{m \text{ times}}),
	\quad
	m\in \IN.
\end{align}

\begin{theorem}[{\cite[Theorem 14.7]{Muj86}}]\label{thm:equivalence_FrechetDer_Hol}
Let $U$ be an open, nonempty subset of $E$. 
For the map {$\mathcal{F}:U \subset E \rightarrow F$} the 
following conditions are equivalent:
\begin{itemize}
\item[(i)] $\mathcal{F}$ is holomorphic.
\item[(ii)] $\mathcal{F}$ is Fr\'echet differentiable.
\item[(iii)] $\mathcal{F}$ is infinitely Fr\'echet differentiable.
\end{itemize}
\end{theorem}

In addition, we recall the following result that will
be of importance of the shape holomorphy analysis of
volumen integral operators.

\begin{lemma}[{\cite[Theorem 1.50]{henriquez2021shape}}]\label{thm:Taylor_holomorphic}
Let $U$ be an open subset of $E$ and let $\mathcal{F}: U \subset E \rightarrow F$
be holomorphic. If the segment joining $r\in U$ and $r+\xi \in U$ is contained
in $U$, then for all $m\in \mathbb{N}_0$ it holds
\begin{align}
	\mathcal{F}(r+\xi)
	=
	\sum_{\ell=0}^{m}
	\frac{1}{\ell!}
	\left(\frac{d^{\ell}}{dr^{\ell}} \mathcal{F}\right)(r,\xi)
	+
	\int\limits_{0}^{1} \frac{(1-\eta)^{m} }{m!} 
        \left(
        		\frac{d^{m+1}}{d r^{m+1}} \mathcal{F}
	\right)(r+\eta \xi,\xi) d\eta.
\end{align}
\end{lemma}

\begin{lemma}[{\cite[Corollary 7.3]{Muj86}}]\label{prop:cauchy_integral_formula}
Let $U$ be an open, nonempty subset of $E$ 
and let {$\mathcal{F}:U \rightarrow F$} be Fr\'echet differentiable. 
Let $r\in U$, $\xi\in E$ and $\vartheta>0$ be such that $r+\sigma \xi \in U$,
for all $\sigma \in \overline{\mathscr{D}}(\vartheta)$. 
Then, for each $m\in \mathbb{N}_0$
holds the Cauchy's integral formula
\begin{align}
	 \frac{d^m}{dr^m}  \mathcal{F}(r,\xi)
	 =
	 \frac{m!}{2\pi \imath} \int\limits_{\snorm{\lambda}=\vartheta} 
	 \frac{\mathcal{F}(r+\lambda \xi)}{\lambda^{m+1}} d\lambda.
\end{align}
\end{lemma}

\begin{lemma}[{\cite[Theorem 3.1.5, item c)]{herve2011analyticity}}]\label{lmm:herve}
Let $U$ be an open, nonempty subset of $E$, and 
let $\left(\mathcal{F}_n\right)_{n\in\mathbb{N}}$ 
be a sequence of holomorphic maps from $U \subset E$ to $F$.
Assume that $\mathcal{F}_n$ converges uniformly to $\mathcal{F}$
in $U$.
Then $\mathcal{F}: U \subset E \rightarrow F$ is holomorphic.
\end{lemma}

%
%

\subsection{Riesz-Thorin Interpolation Theorem}
\label{sec:rt_theorem}
Let $\D$ be a bounded Lipschitz domain in $\IR^d$, $d \in \{2,3\}$.
We consider a volume integral operator of the form
\begin{equation}
	(\OT \varphi)(\bx)
	\coloneqq
	\int\limits_{\D}
	{\sf t}(\bx,\by)\varphi(\by)\mathrm{d}\by,
	\quad
	\bx
	\in
	\D,
\end{equation}
an integral operator defined in $\D$ with 
$\mathsf{t}: \left(\D \times \D \right)^\star \rightarrow \IC$.
To show that $\OT\colon L^2(\D)\to L^2(\D)$
defines a bounded linear operator
we rely on the \emph{Riesz-Thorin interpolation theorem}
\cite[Theorem 2.b.14]{LT77}:
The operator norm of $\OT: L^2(\D)
\rightarrow  L^2(\D)$
is bounded according to
\begin{align}\label{eq:rieszthorin}
	\norm{
		\OT
	}_{\mathscr{L}\left(L^2(\D), L^2(\D)\right)}
	\leq
	\norm{
		\OT
	}_{\mathscr{L}\left(L^1(\D), L^1(\D)\right)}^{\half}
	\norm{
		\OT
	}_{\mathscr{L}\left(L^\infty(\D),L^\infty(\D)\right)}^{\half},
\end{align}
and we estimate the right-hand side of \cref{eq:rieszthorin} as follows
\begin{align}
	\norm{
		\OT
	}_{\mathscr{L}\left(L^1(\D), L^1(\D)\right)}
	&=
	\underset{\by \in {\D}}{\operatorname{ess} \sup} 
	\int\limits_{\D}|\mathsf{t}(\bx,\by)| \mathrm{d}\bx,
	\label{eq:rt_bounds_1}
	\\
	\norm{
		\OT
	}_{\mathscr{L}\left(L^\infty(\D), L^\infty(\D)\right)}
	&=
	\underset{\bx \in {\D}}{\operatorname{ess \sup}}
	\int\limits_{\D}|\mathsf{t}(\bx,\by)| \mathrm{d}\by.
	\label{eq:rt_bounds_2}
\end{align}

\section{Volumen Integral Equations}
\label{sec:vol_int_eq}
The model problem to be considered in this work 
the Helmholtz equation $\IR^d, d = 2,3$ with variable
coefficients in a compact subset of the entire space. 
More precisely, we seek $u:\IR^d \rightarrow \IC$
satisfying
\begin{equation}\label{eq:Helmholtz}
	 -\Delta u(\bx) - \kappa(\bx)^2u(\bx) 
	 = 
	 f(\bx), 
	 \quad 
	 \bx \in \IR^d,
\end{equation}
where $\kappa \in L^{\infty}(\IR^d)$ is 
such that $\kappa(\bx) \equiv \kappa_0 > 0 $ for a.e. $\bx \in \IR^d \setminus \overline{\text{D}}$,
with $\text{D}$ a bounded subset of $\IR^d$.
In addition, we assume $f\in L^2_{\text{comp}}(\IR^d)$ to be compactly supported in $\text{D}.$

As it is customary in the treatment of boundary value
problems in unbounded domains, we decompose the
total field $u$ as $u= u^s + \uinc$, where $u^s$ satisfies
the Sommerfeld radiation condition
\begin{equation}\label{eq:radiation_conditions}
    \lim_{r\rightarrow \infty}r^{\frac{d-1}{2}} \left(\dfrac{\partial u^s}{\partial r} - \imath\kappa_0u^s \right) = 0, \quad r = \norm{\bx},
\end{equation}
and the incident field $\uinc$ satisfies the homogeneous
Helmholtz equation, i.e.
\begin{equation}\label{eq:uincHelmholtz}
    -\Delta \uinc - \kappa_0^2\uinc = 0 \quad \text{ in }\IR^d.
\end{equation}

Set $\beta(\bx) := \kappa(\bx)^2 - \kappa_0^2$ for $\bx \in \IR^d$.
Consequently, 
\cref{eq:Helmholtz} yields
\begin{equation}\label{eq:new_Helmholtz}
	-\Delta u(\bx) 
	- 
	\kappa_0 ^2u(\bx) 
	= 
	\beta(\bx)u(\bx) 
	+ 
	f(\bx),
	\quad 
	\bx
	\in
	\IR^d.
\end{equation}
Observe that $\beta$ is a compactly supported function in $\overline{\text{D}}$.

Given $\omega>0$, let $G^{(\omega)}: \IR^d \backslash \{{\bf 0}\}\rightarrow \IC$
be the fundamental solution of the Helmholtz equation with
wavenumber $\omega > 0$ in $\IR^d$
\begin{equation}\label{eq:green_G_omega}
	G^{(\omega)}(\bz)
	\coloneqq
	\left \{
	\begin{array}{cl}
		\frac{\imath}{4} \text{H}^{(1)}_0(\omega \norm{\bz}) & d=2 \\
		\frac{1}{4\pi} \frac{\exp(\imath \omega \norm{\bz})}{\norm{\bz}} & d=3 \\
	\end{array}
	\right.
	,
	\quad
	\bz
	\in
	\IR^d \backslash \{{\bf 0}\},
\end{equation}
where $ \text{H}^{(1)}_0$ corresponds to the Hankel function of the first kind and order zero.

The Newton potential is defined for $v\in C^{\infty}_{\text{comp}}(\IR^d)$ as
\begin{equation}
	\left(
		\N^{(\omega)}v
	\right)(\bx)
	\coloneqq
	\int\limits_{\IR^d}
	G^{(\omega)}
	(\bx-\by)v(\by)\text{d}\by,
	\quad
	\bx \in \mathbb{R}^d,
\end{equation}
where $\omega \in \mathbb{R}$.
In particular, it can be extended as a continuous operator 
\begin{equation}
	\N^{(\omega)}: 
	H^s_{\text{comp}}(\IR^d) 
	\rightarrow 
	H^{s+2}_{\text{loc}}(\IR^d), \quad \text{for all }s\in\IR.
\end{equation}
For the exact definition of these spaces, we refer to \cite[Definition 2.6.1 and Definition 2.6.5]{SS10}
As a consequence of Green's identity \cite[Theorem~8.3]{CK12},
a solution of \cref{eq:radiation_conditions}--\cref{eq:new_Helmholtz}
admits the following representation
\begin{equation}\label{eq:LS-rep}
	u 
	= 
	\uinc 
	+ 
	\N^{(\kappa_0)}(\beta u) 
	+ 
	\N^{(\kappa_0)}
	f, 
	\quad 
	\text{ in }\IR^d.
\end{equation}
Note that \cref{eq:LS-rep} is both a representation formula
and an integral equation.
Restricting the problem to the domain $\text{D}$ and
assuming $f\equiv 0,$ we obtain the so-called
Lippmann-Schwinger equation or 
Helmholtz VIE:
Find $u \in H^1(\text{D})$ satisfying
\begin{equation}\label{eq:LS}
	u - \mathsf{A}u = \uinc
	\quad
	\text{ in }
	\text{D}.
\end{equation}
where 
\begin{equation}\label{eq:VIO}
	\mathsf{A}:
	\N^{(\kappa_0)}(\beta \cdot)|_{\text{D}}: \ H^1(\text{D}) \rightarrow H^1(\text{D})
\end{equation}
is a compact operator. \Cref{eq:LS} is therefore studied in $H^1(\text{D})$. 
The variational formulation of \cref{eq:LS} reads as follows.

\begin{problem}[Variational Formulation of the Lippmann-Schwinger Equation]\label{LSvar}
Let $\normalfont\text{D}\subset \IR^d$ be a bounded Lipschitz domain. 
Let $\kappa\in L^{\infty}(\IR^d), \ \kappa_0 > 0$ 
such that $\normalfont\text{supp}(\kappa-\kappa_0)\subset \text{D}$,
and let $\mathsf{A}$ be as in \cref{eq:VIO}.
Given an incident field $\uinc$ satisfying \cref{eq:uincHelmholtz}, 
we seek $u\in H^1(\text{D})$ such that
\begin{equation}\label{eq:bilinearLS}
	\mathsf{a}(u, v) 
	\coloneqq \dual{u}{v}_{\normalfont\text{D}} - \dual{\mathsf{A}u}{v}_{\normalfont\text{D}}
	= \dual{\normalfont\uinc}{v}_{\normalfont\text{D}},
\end{equation}
holds for all $v\in \widetilde{H}^{-1}(\normalfont\text{D}).$
\end{problem}
In this particular setting, \Cref{LSvar} is well-posed.
Indeed, the bilinear form $\mathsf{a}$ in \cref{eq:bilinearLS}
is bounded and satisfies the T-coercivity property (see \cite[Theorem~1]{ciarlet2012t}, \cite[Lemma~1.4.3]{labarca2024coupled})
Moreover, the formulation in \cref{eq:LS} is equivalent
to \cref{eq:Helmholtz}--\cref{eq:radiation_conditions}.

\subsection{Numerical Discretization}\label{sec:num_dict}
Let $ \{\mathcal{T}_h\}_{0<h<h_0} $ be a sequence of globally quasi-uniform, shape-regular family of triangular meshes of $ \text{D}$. We choose finite element spaces $ V_h := V_h(\mathcal{T}_h) \subset H^1(\text{D}) $ of piecewise linear functions on $ \mathcal{T}_h $. We also use the same finite dimensional space $ V_h $ as a conforming subspace of $ \widetilde{H}^{-1}(\text{D}). $ 

We proceed to address the stability of discrete counterpart of \cref{LSvar} in 
the finite dimensional space $V_h$.
To this end, we establish a \emph{discrete} inf-sup condition for
the bilinear form $\mathsf{a}: H^1(\text{D}) \times \widetilde{H}^{-1}(\text{D}) \rightarrow \IC$.
As discussed in \cite[Section~2.1]{steinbach2003stability}, 
there exists an $L^2(\text{D})$-orthogonal projection operator 
$ Q_h: H^1(\text{D})\rightarrow V_h \subset H^1(\text{D})$ such that
\begin{equation}
	\langle Q_hu, w_h\rangle_{\text{D}} 
	= 
	\langle u, w_h\rangle_{\text{D}}, \quad \text{for all }w_h\in V_h,
\end{equation}
and satisfying the following properties:
For all $u\in H^1(\text{D})$ it holds
\begin{align}
	\norm{Q_h u}_{H^1(\text{D})} 
	\leq 
	c_S \norm{u}_{H^1(\text{D})}
	\quad
	\text{and}
	\quad
	\norm{u - Q_hu}_{L^2(\text{D})} 
	\leq 
	c_1 h|u|_{H^1(\text{D})},\label{eq:Q_error}
\end{align}
where $ c_S, c_1>0 $ depend only on the shape-regularity and quasi-uniformity measure of $ \mathcal{T}_h, $ but not on the parameter $ h. $ This is equivalent to the following result

\begin{proposition}[Discrete inf-sup condition for a]\label{th:discrete-inf-sup}
Consider the setting of \cref{LSvar}.
There exists $c_{Q}>0 $ such that for all $h>0$ it holds
\begin{equation}
	c_{Q} 
	\leq 
	\inf\limits_{0\neq u_h\in V_h}
	\sup\limits_{0\neq v_h\in V_h}
	\dfrac{
		\Re
		\left\{
			\langle v_h, u_h 
			\rangle_{\normalfont\text{D}}
		\right\}
	}{
		\norm{u_h}_{H^1(\normalfont\text{D})}
		\norm{v_h}_{\widetilde{H}^{-1}(\normalfont\text{D})}
	}.
\end{equation}
\end{proposition}

\cref{th:discrete-inf-sup} entails that, up to a compact perturbation,
the bilinear form $\mathsf{a}(\cdot, \cdot) : V_h \times V_h \rightarrow \mathbb{C}$
satisfies a discrete inf-sup condition.
Now, we are in position to establish the main result 
concerning the Galerkin discretization of \cref{LSvar}.

\begin{theorem}\label{th:galerkin_qo}
There is $ h_0 >0 $ and a constant $c_{\textsf{qo}} >0$ 
independent of $ h $ such that for each $0<h<h_0$ there exists a unique
$u_h \in V_h$ solution of $\mathsf{a}(u_h, v_h) = \dual{\normalfont\uinc}{v_h}_{\normalfont\text{D}}$, for all 
$v_h \in V_h$, and satisfying
\begin{equation}
	\norm{u - u_h}_{H^1(\normalfont\text{D})}
	\leq 
	c_{\textsf{qo}} 
	\inf\limits_{w_h \in V_h}
	\norm{u - w_h}_{H^1(\normalfont\text{D})},
\end{equation}
where $u \in H^1(\normalfont\text{D})$ is the solution of \cref{LSvar}.
\end{theorem}

\section{Shape Holomorphy}
\label{sec:shape_holomorphy}
In this section, we study the holomorphic dependence of a collection of 
volume integral operators, as the one introduced in \Cref{sec:vol_int_eq}, upon a
suitable family of domain transformations.

\subsection{Holomorphic Volume Integral Operators}
\label{sec:holomorphic_volume_operators}
Let $\D$ be simply connected, bounded Lipschitz domain in 
$\IR^{d}$, $d\in \{2,3\}$. For $\varphi \in L^2(\D)$ set
\begin{align}\label{eq:volume_operator}
	\left(
		\mathsf{P}_{\br}
		\varphi
	\right)
	(\bx)
	=
	\int\limits_{\D}
	\mathsf{p}_{\br}(\bx,\by)
	\varphi(\by)
	\dd \by,
	\quad
	\bx
	\in 
	\D,
\end{align}
where we have assumed that $\mathsf{p}_{\br}$ depends
on the parameter $\br \in \mathfrak{T}$, being $\mathfrak{T}$ a compact 
subset of a \emph{complex} Banach space $X$. 

In the following result we address the holomorphic
dependence of the volume integral operator $\mathsf{P}_{\br}$
(as an element of suitable space of bounded linear operators)
upon  $\br \in \mathfrak{T}$.
To this end, we first address the case when $\mathsf{p}_{\br}$
is bounded in $\D$. This result serves as a fundamental building
block to establish shape holomorphy of the VIO introduced in \Cref{sec:vol_int_eq}.

Let $\mathfrak{T}$ be a compact subset of $X$.
Given $\varepsilon>0$ set
\begin{equation}
	\mathfrak{T}_{\varepsilon}
	\coloneqq
	\left\{
		\br 
		\in 
		X
		:
		\exists\,
		\widetilde{\br} \in \mathfrak{T}
		\text{ such that }
		\norm{\br - \widetilde{\br}}_{X}
		<
		\varepsilon
	\right\}
\end{equation}

\begin{theorem}\label{eq:holomorphy_bounded_kernel}
Let $\D$ be simply connected, bounded Lipschitz domain in 
$\IR^{d}$, $d\in \{2,3\}$. 
Assume that there exists $\delta>0$ such that the following
conditions are satisfied:
\begin{itemize}
	\item[(i)]
	For each $\br \in \mathfrak{T}$ the map 
	$\mathfrak{T} \ni \br \mapsto \mathsf{p}_{\br} \in L^{\infty}(\D \times \D)$
	admits a uniformly bounded holomorphic extension into $\mathfrak{T}_\delta$
	denoted by
	$\mathfrak{T}_\delta \ni \br \mapsto \mathsf{p}_{\br,\IC} \in L^{\infty}(\D \times \D)$.
	\item[(ii)]
	The extension of 
	$\mathsf{P}_{\br}: L^2(\D) \rightarrow L^2(\D)$
	to $\mathfrak{T}_\delta$
	defined as
	\begin{align}\label{eq:extension_P}
		\left(
			\mathsf{P}_{\br,\IC}
			\,
			\varphi
		\right)
		(\bx)
		=
		\int\limits_{\D}
		\mathsf{p}_{\br,\IC}(\bx,\by)
		\varphi(\by)
		\dd \by,
		\quad
		\bx
		\in 
		\D
	\end{align}
	is uniformly bounded upon $\mathfrak{T}_\delta$,
	i.e. there exists $C_{\mathsf{P}}>0$ such that
	\begin{equation}
		\sup_{\br \in \mathfrak{T}_\delta}
		\norm{
			\mathsf{P}_{\br,\IC}
		}_{\mathscr{L}(L^2(\D),L^2(\D))}
		\leq
		C_{\mathsf{P}}
	\end{equation}
\end{itemize}
Then, the map
\begin{equation}
	\mathcal{P}:
	\mathfrak{T}
	\rightarrow
	\mathscr{L}
	\left(
		L^2(\D)
		,
		L^2(\D)
	\right):
	\br
	\mapsto
	\mathsf{P}_{\br}
\end{equation}
admits a bounded holomorphic extension into $\mathfrak{T}_\delta$.
\end{theorem}

\begin{proof}
This results follows from the exact same arguments used in the proofs of
\cite[Theorem 3.12]{HS21}, \cite[Theorem 3.1]{DH23}, and \cite[Theorem 4.5]{PHJ23}.
For the sake of brevity, we skip it.
\end{proof}

The previous result addresses the holomorphic dependence of the 
operator provided that for each $\br \in \mathfrak{T}_\delta$ 
we have $\mathsf{p}_\br \in L^\infty(\D\times \D)$.
Let us set
\begin{equation}
	({\D} \times {\D})^\star 
	\coloneqq 
	\left\{
		(\bx,\by) 
		\in 
		{\D} \times {\D}: 
		\bx\neq\by 
	\right\}.
\end{equation}
Next, we proceed to extend this results to
$\mathsf{p}_\br: ({\D} \times {\D})^\star \rightarrow \IC$,
which possibly has a singularity at $\bx = \by$.
A key instrument to prove this result is \Cref{lmm:herve}

\begin{theorem}\label{eq:holomorphic_L2}
Let $\D$ be simply connected, bounded Lipschitz domain in 
$\IR^{d}$, $d\in \{2,3\}$. 
Assume that there exists $\delta>0$ such that the following
conditions are satisfied:
\begin{itemize}
	\item[(i)]
	For each $(\bx,\by) \in ({\D} \times {\D})^\star $
	the map 
	\begin{equation}
		\mathfrak{T} \ni \br \mapsto \mathsf{p}_{\br}(\bx,\by) \in \IC
	\end{equation}
	admits a uniformly bounded holomorphic extension into $\mathfrak{T}_\delta$
	denoted by
	$\mathfrak{T}_\delta \ni \br \mapsto \mathsf{p}_{\br,\IC}(\bx,\by) \in \IC$.
	\item[(ii)]
	There exists $C(\delta,\mathfrak{T})>0$ and $\nu \in[0,d)$, depending on 
	$\delta>0$ and $\mathfrak{T}$ only,
	such that for each $\br\in \mathfrak{T}_\delta$ it holds
	\begin{equation}\label{eq:singular_kernel}
		\snorm{
			\mathsf{p}_{\br,\IC}
			(\bx,\by)
		}
		\leq
		\frac{
			C(\delta,\mathfrak{T})
		}{
			\norm{\bx-\by}^\nu
		},
		\quad
		(\bx,\by) \in ({\D} \times {\D})^\star.
	\end{equation}
\end{itemize}
Then, the map
\begin{equation}\label{eq:P_shape_hol_Op}
	\mathcal{P}:
	\mathfrak{T}
	\rightarrow
	\mathscr{L}
	\left(
		L^2(\D)
		,
		L^2(\D)
	\right):
	\br
	\mapsto
	\mathsf{P}_{\br}
\end{equation}
admits a bounded holomorphic extension into $\mathfrak{T}_\delta$.
\end{theorem}

\begin{proof}
We take our cue from \cite[Theorem 3.2]{DH23}.
Let $\chi \in \mathscr{C}^\infty([0,\infty))$ be a smooth function
satisfying the following properties: $\chi(t) =0$ for $t\in[0,\half]$;
$\chi(t) =1$ for $t\geq 1$; and $\chi(t) \in [0,1]$ for $t\in [0,\infty)$.

For each $n\in \IN$, we set
\begin{align}\label{eq:P_n}
	\left(
		{\mathsf{P}}^{(n)}_{\br,\IC}
		\,
		\varphi
	\right)
	({\bf x})
	\coloneqq
	\int\limits_{\D}
	\mathsf{p}_{\br,\IC}(\bx,\by)
	\chi
	\left(
		n\norm{\bx-\by}^\nu
	\right)
	\varphi(\by)
	\dd{\by},
	\quad
	{\bx}
	\in 
	{\D}.
\end{align}
Set
\begin{equation}
	\mathsf{p}^{(n)}_{\br,\IC}(\bx,\by)
	\coloneqq
	\mathsf{p}_{\br,\IC}(\bx,\by)
	\chi
	\left(
		n\norm{\bx-\by}^\nu
	\right),
	\quad
	(\bx,\by) \in {\D} \times {\D}.
\end{equation}
Clearly for each $\br \in \mathfrak{T}_\delta$ one has
$\mathsf{p}^{(n)}_{\br,\IC} \in L^\infty({\D} \times {\D})$, thus
yielding ${\mathsf{P}}^{(n)}_{\br,\IC} \in \mathscr{L}(L^2(\D),L^2(\D))$
for each $n\in \IN$. 
Also, due to \cref{eq:singular_kernel},
\cref{eq:rt_bounds_1}, and \cref{eq:rt_bounds_2}
one obtains
\begin{equation}
\begin{aligned}
	\norm{
		{\mathsf{P}}_{\br,\IC}
	}_{\mathscr{L}\left(L^1(\D), L^1(\D)\right)}
	&
	=
	\underset{\by \in {\D}}{\operatorname{ess} \sup} 
	\int\limits_{\D}\snorm{
		\mathsf{p}_{\br,\IC}(\bx,\by)
	}\mathrm{d}\bx
	\leq
	C(\delta,\mathfrak{T})
	\;
	\underset{\by \in {\D}}{\operatorname{ess} \sup} 
	\int\limits_{\D}
	\frac{
		\mathrm{d}\bx
	}{
		\norm{\bx -\by}^\nu
	}
	<
	\infty,
\end{aligned}
\end{equation}
and
\begin{equation}
\begin{aligned}
	\norm{
		{\mathsf{P}}_{\br,\IC}
	}_{\mathscr{L}\left(L^\infty(\D), L^\infty(\D)\right)}
	&
	=
	\underset{\bx \in {\D}}{\operatorname{ess} \sup} 
	\int\limits_{\D}\snorm{
		\mathsf{p}_{\br,\IC}(\bx,\by)
	}\mathrm{d}\by
	\leq
	C(\delta,\mathfrak{T})
	\;
	\underset{\bx \in {\D}}{\operatorname{ess} \sup} 
	\int\limits_{\D}
	\frac{
		\mathrm{d}\bx
	}{
		\norm{\bx -\by}^\nu
	}
	<\infty.
\end{aligned}
\end{equation}
Thus, $\mathsf{P}_{\br,\IC} \in \mathscr{L}(L^2(\D),L^2(\D))$,
and 
\begin{equation}
	\sup_{\br \in \mathfrak{T}_\delta}
	\norm{
		{\mathsf{P}}_{\br,\IC}
	}_{\mathscr{L}\left(L^2(\D), L^2(\D)\right)}
	\leq
	\widetilde{C}(\delta,\mathfrak{T}).
\end{equation}
Next, we show that
\begin{equation}\label{eq:unform_convergence_P}
	\limsup_{n\rightarrow \infty}
	\sup_{\br\in \mathfrak{T}_\delta}
	\norm{
		{\mathsf{P}}^{(n)}_{\br,\IC}
		-
		{\mathsf{P}}_{\br,\IC}
	}_{\mathscr{L}(L^2(\D),L^2(\D))}
	=
	0.
\end{equation}
For $\varphi \in L^2(\D)$, $\bx \in \D$,
and $n\in \IN$ one has
\begin{align}\label{eq:difference_operators}
	\left(
		{\mathsf{P}}_{\br,\IC}
		\,
		\varphi
	\right)
	(\bx)
	-
	\left(
		{\mathsf{P}}^{(n)}_{\br,\IC}
		\,
		\varphi
	\right)
	(\bx)
	=
	\int\limits_{\D}
	\left(
		1
		-
		\chi
		\left(
			n\norm{\bx-\by}^\nu
		\right)
	\right)
	\mathsf{p}_{\br,\IC}(\bx,\by)
	\varphi(\by)
	\dd \by.
\end{align}
The right-hand side of \cref{eq:difference_operators}
is itself a volumen integral operator. 
Consequently, for $n\in \IN$,
for each $r\in \mathfrak{T}_\delta$
and recalling the Riesz-Thorin interpolation theorem of 
\Cref{sec:rt_theorem} we get
\begin{equation}
\begin{aligned}
	\norm{
		{\mathsf{P}}_{\br,\IC}
		-
		{\mathsf{P}}^{(n)}_{\br,\IC}
	}_{\mathscr{L}\left(L^1(\D), L^1(\D)\right)}
	&
	=
	\underset{\by \in {\D}}{\operatorname{ess} \sup} 
	\int\limits_{\D}
	\snorm{
		1
		-
		\chi
		\left(
			n\norm{\bx-\by}^\nu
		\right)
	}
	\snorm{
		\mathsf{p}_{\br,\IC}(\bx,\by)
	}
	\dd \bx
	\\
	&
	\leq
	C(\delta,\mathfrak{T})
	\;
	\underset{\by \in {\D}}{\operatorname{ess} \sup} 
	\int\limits_{\D  \cap B(\by,n^{-\frac{1}{\nu}})}
	\frac{
		\mathrm{d}\bx
	}{
		\norm{\bx -\by}^\nu
	}
	\dd \bx
	\\
	&
	\leq
	C(\delta,\mathfrak{T})
	\omega_d
	\;
	\int\limits_{0}^{n^{-\frac{1}{\nu}}}
	\varrho^{-\nu+d-1}
	\dd \varrho
	\\
	&
	\leq
	C(\delta,\mathfrak{T})
	\omega_d
	\;
	\left.
		\frac{
			\varrho^{-\nu+d}
		}{
			-\nu+d
		}
	\right\rvert_{0}^{n^{-\frac{1}{\nu}}}
	=
	\frac{
		C(\delta,\mathfrak{T})
	}{
		-\nu+d
	}
	\omega_d
	n^{\frac{\nu-d}{\nu}}
\end{aligned}
\end{equation}
where we have used item (ii) in the assumptions,
and the fact that $\nu \in [0,d)$, and $\omega_d$
is the surface of the unit sphere in $\IR^d$.
Using the this tools, one can prove 
\begin{equation}
		\norm{
		{\mathsf{P}}_{\br,\IC}
		-
		{\mathsf{P}}^{(n)}_{\br,\IC}
	}_{\mathscr{L}\left(L^\infty(\D), L^\infty(\D)\right)}
	\leq
	\frac{
		C(\delta,\mathfrak{T})
	}{
		-\nu+d
	}
	\omega_d
	n^{\frac{\nu-d}{\nu}}.
\end{equation}
Observe that these bounds are uniform on $\mathfrak{T}_\delta$.
Consequently, by recalling again that $\nu \in [0,d)$ we can show
that \cref{eq:unform_convergence_P} holds.
Observe that for each $n\in \IN$
and each $\br \in \mathfrak{T}_\delta$
we have that $\mathsf{p}^{(n)}_{\br,\IC} \in 
L^\infty(\D \times \D)$ with $\delta>0$
independent of $n\in \IN$.

We proceed to prove that for each finite $n\in \IN$
the map
\begin{align}\label{eq:holomorphic_kernel}
	\mathfrak{T}_\delta
	\ni
	\br
	\mapsto
	\mathsf{p}^{(n)}_{\br,\IC}
	\in 
	L^\infty(\D \times \D)
\end{align}
is holomorphic and uniformly bounded
with $\delta>0$ as in \cref{eq:prop_existence_delta}.

Let $\varepsilon\in(0,\delta)$. 
Using \cref{thm:Taylor_holomorphic,prop:cauchy_integral_formula},
for each $r\in \mathfrak{T}_\delta$
and for each fixed
$(\bx,\by) \in (\D \times \D)^\star$
it holds
\begin{align}\label{eq:taylor_expansio_kernel}
	{\mathsf{p}}_{\br+\bxi,\IC}(\bx,\by) 
	=
	{\mathsf{p}}_{\br,\IC}(\bx,\by) 
	+
	\left(
		\frac{{d}}{{d} \br}{\mathsf{p}}_{\cdot,\IC}
	\right)(\br,\bxi)(\bx,\by) 
	+
	\frac{1}{\pi \imath} 
	\int\limits_{\snorm{\lambda} 
	= 
	\vartheta_\bxi} 
	\frac{{\mathsf{p}}_{\br+\lambda\bxi,\IC}(\bx,\by)}{\lambda^{3}} 
	\text{d}\lambda,
\end{align}
where $0 \neq \bxi \in \mathfrak{T}_\delta$
satisfying
$\norm{\bxi}_{X}<\delta-\varepsilon$
and
\begin{align}\label{eq:vartheta}
	0<
	\vartheta_\bxi 
	\coloneqq
	\frac{\delta-\varepsilon}{\norm{\bxi}_{X}} 
	- 
	1.
\end{align}
Multiplying \cref{eq:taylor_expansio_kernel}
by $\chi(n\norm{\bx-\by}^\nu)$
we obtain that for each $r\in \mathfrak{T}_\delta$,
$n\in \IN$, and $(\bx,\by) \in \D \times \D$
it holds
\begin{equation}\label{eq:taylor_kernel_p_n}
\begin{aligned}
	\left(
	{\mathsf{p}}_{\br+\bxi,\IC}(\bx,\by) 
	-
	{\mathsf{p}}_{\br,\IC}(\bx,\by) 
	-
	\left(
		\frac{d}{d\br}{\mathsf{p}}_{\cdot,\IC}
	\right)(\br,\bxi)(\bx,\by) 
	\right)
	&
	\chi
	\left(
		n\norm{\bx-\by}^\nu
	\right)  \\
	&
	\hspace{-2.5cm}
	= 
	\frac{1}{\pi \imath} 
	\int\limits_{\snorm{\lambda} 
	= 
	\vartheta_\bxi}
	\frac{
		{\mathsf{p}}_{\br+\lambda\bxi,\IC}(\bx,\by)
	\chi
	\left(
		n\norm{\bx-\by}^\nu
	\right)
		}{\lambda^{3}} d\lambda \\
	&
	\hspace{-2.5cm}
	= 
	\frac{1}{\pi \imath} 
	\int\limits_{\snorm{\lambda} 
	=
	\vartheta_\bxi}
	\frac{
		n 
		\norm{\bx-\by}^\nu
		\mathsf{p}_{\br+\lambda\bxi,\IC}(\bx,\by)
		f^{(n)}(\bx,\by)
	}{\lambda^{3}} d\lambda,
\end{aligned}
\end{equation}
where, for each $n\in \IN$, 
$f^{(n)}: \D \times \D \rightarrow \IR$ is defined as
\begin{equation}
	f^{(n)}
	(\bx,\by) 
	\coloneqq
	\left\{
	\begin{array}{cl}
	\frac{
	\chi
	\left(
		n\norm{\bx-\by}^\nu
	\right)-\chi(0)}{n\norm{\bx-\by}^\nu},&
	(\bx,\by)  
	\in 
	({\D} \times {\D})^\star, \\
	\chi'(0),&
	(\bx,\by)  
	\notin 
	(\D \times \D)^\star,
	\end{array}
	\right.
\end{equation}
and 
\begin{align}
	\norm{f^{(n)}}_{L^\infty(\D \times \D)}
	\leq
	\sup_{t\in[0,\infty)}
	\snorm{\chi'(t)}
	<
	\infty.
\end{align}
Then, for each $n\in \IN$
\begin{align}
	\norm{
	{\mathsf{p}}^{(n)}_{\br+\bxi,\IC}
	-
	{\mathsf{p}}^{(n)}_{\br,\IC} 
	-
	\left(
		\frac{d}{d\br}{\mathsf{p}^{(n)}}_{\cdot,\IC}
	\right)(\br,\bxi)
	}_{L^\infty({\D}\times {\D})}
	\leq
	2
	\frac{n}{ \vartheta^2_\bxi}
	\sup_{t\in[0,\infty)}
	\snorm{\chi'(t)}
	C(\delta,\mathfrak{T}),
\end{align}
with $C(\delta,\mathfrak{T})>0$ as in \cref{eq:bound_kernel}
and with 
\begin{align}
	\left(
		\frac{d}{d\br}{\mathsf{p}^{(n)}_{\cdot,\IC}}
	\right)(\br,\bxi)
	\left(
		{\bx},{\by}
	\right)
	=
	\left(
		\frac{d}{d\br}{\mathsf{p}_{\cdot,\IC}}
	\right)(\br,\bxi)
	\left(
		{\bx},{\by}
	\right)
	\chi(n\norm{\bx-\by}^\nu),
	\quad
	(\bx,\by)
	\in
	\D \times \D.
\end{align}
Observing that
\begin{align}
	\frac{1}{\vartheta^2_\bxi} 
	= 
	o
	\left(
		\norm{\bxi}_{X}
	\right)
\end{align}
we get the desired result. 

We remark that this result is 
valid for each finite $n\in \IN$. As we let $n\rightarrow \infty$
the holomorphy property does not hold anymore due to the singular
behaviour at $\bx \neq \by$.

Now, we claim that for each $n\in \IN$
the map 
\begin{align}
	\mathcal{P}^{(n)}:
	\mathfrak{T}_\delta
	\rightarrow
	\mathscr{L}
	\left(	
		L^2(\text{D})
		,
		L^2(\text{D})
	\right)
	:
	\br
	\mapsto
	{\mathsf{P}}^{(n)}_{\br,\IC}
\end{align}
is holomorphic with Fr\'echet derivative at $\br \in \mathfrak{T}_\delta$
along the direction $\bxi \in X$ given by
\begin{equation}
	\left(
		\frac{d}{d\br}
		\mathcal{P}^{(n)}
		\varphi
	\right)
	(\br,\bxi)
	(\bx)
	=
	\int\limits_{\D}
	\left(
		\frac{d}{d\br}{\mathsf{p}^{(n)}_{\cdot,\IC}}
	\right)(\br,\bxi)
	\left(
		{\bx},{\by}
	\right)
	\varphi(\by)
	\dd{\by},
	\quad
	{\bx}
	\in 
	{\D}.
\end{equation}
To this end, we invoke \cref{eq:holomorphy_bounded_kernel}.
Consequently, it remains to verify the assumptions of said result. 
Item (i) has been verified since for each $\br \in \mathfrak{T}_\delta$
and for each $n\in \mathbb{N}$ we have that
$\mathsf{p}^{(n)}_{\br,\IC} \in L^\infty(\D\times\D)$.
Also, we have shown here that the map in \cref{eq:holomorphic_kernel}
is holomorphic. The uniform boundedness in item (ii) has been
assumed in this result as well.

Finally, it follows from \Cref{lmm:herve} that the map introduced in \cref{eq:P_shape_hol_Op}
admits a bounded holomorphic extension into $\mathfrak{T}_\delta$.
\end{proof}

However, since we are interested in studying
as well the VIOs introduced in
\Cref{sec:vol_int_eq} as elements of the complex Banach space
$\mathscr{L}(L^2(\text{D}),H^1(\text{D}))$ we introduce 
the following result.

\begin{theorem}\label{eq:holomorphic_H1}
Consider the assumptions of \cref{eq:holomorphic_L2}, however with $\nu \in [0,d-1)$.
In addition, we assume the following.
\begin{itemize}
\item[(i)]
For each $(\bx,\by) \in ({\D} \times {\D})^\star $
the map 
\begin{equation}
	\mathfrak{T} \ni \br \mapsto \nabla_\bx \mathsf{p}_{\br}(\bx,\by) \in \IC
\end{equation}
admits a uniformly bounded holomorphic extension into $\mathfrak{T}_\delta$
denoted by
\begin{equation}
	\mathfrak{T}_\delta \ni \br \mapsto \nabla_{\bx} \mathsf{p}_{\br,\IC}(\bx,\by) \in \IC.
\end{equation}
\item[(ii)]
For $\zeta \in[0,d)$ and $\delta>0$
as in \cref{eq:holomorphic_L2} we have that for 
a.e. each $\br \in \mathfrak{T}_\delta$ it holds
\begin{equation}\label{eq:bound_kernel}
	\norm{
		\nabla_{\bx}
		\mathsf{p}_{\br,\IC}
			(\bx,\by)
	}
	\leq
	\frac{
		\check{C}(\delta,\mathfrak{T})
	}{
		\norm{\bx-\by}^{\zeta}
	},
	\quad
	\text{for a.e. }
	(\bx,\by) \in ({\D} \times {\D})^\star,
\end{equation}
where $\nabla_{\bx}$ signifies the weak gradient with 
respect to $\bx$.
\end{itemize}

Then, there exists $\delta>0$ such that the map
\begin{equation}
	\mathcal{P}:
	\mathfrak{T}
	\rightarrow
	\mathscr{L}
	\left(
		L^2(\D)
		,
		H^1(\D)
	\right):
	\br
	\mapsto
	\mathsf{P}_{\br}.
\end{equation}
admits a bounded holomorphic extension into $\mathfrak{T}_\delta$.
\end{theorem}

\begin{proof}
For each $n\in \IN$, $\br \in \mathfrak{T}$, and any $\varphi \in L^2(\D)$
\begin{equation}\label{eq:der_P_n}
\begin{aligned}
	\nabla_{\bx}
	\left(
		{\mathsf{P}}^{(n)}_{\br,\IC}
		\,
		\varphi
	\right)
	({\bf x})
	=
	&
	\int\limits_{\D}
	\left(
		\nabla_{\bx}
		\mathsf{p}_{\br,\IC}(\bx,\by)
		\chi
		\left(
			n\norm{\bx-\by}^\nu
		\right)
	\right.
	\\
	&
	\left.
		+
		n \nu
		\mathsf{p}_{\br,\IC}(\bx,\by)
		\chi'
		\left(
			n\norm{\bx-\by}^\nu
		\right)
		\norm{\bx-\by}^{\nu-2}
		(\bx-\by)
	\right)
	\varphi(\by)
	\dd{\by},
	\quad
	{\bx}
	\in 
	{\D},
\end{aligned}
\end{equation}
with ${\mathsf{P}}^{(n)}_{\br,\IC}:L^2(\text{D}) \rightarrow L^2(\text{D})$ as in \cref{eq:P_n}.

Due to the properties of the function $\chi \in \mathscr{C}^\infty([0,\infty))$
we have that for each $n \in \IN$ and for each $\br \in \mathfrak{T}_\delta$ the operator
${\mathsf{P}}^{(n)}_{\br,\IC}:L^2(\D) \rightarrow H^1(\D)$ is linear
and bounded.
Set 
\begin{equation}
	\left(
		\mathsf{G}_{\br,\IC}
		\,
		\varphi
	\right)
	(\bx)
	=
	\int\limits_{\D}
	\nabla_{\bx}
	\mathsf{p}_{\br,\IC}(\bx,\by)
	\varphi(\by)
	\dd{\by},
	\quad
	{\bx}
	\in 
	{\D}.	
\end{equation}
As a consequence of \cref{eq:bound_kernel} and the Riesz-Thorin theorem stated in \Cref{sec:rt_theorem}, for each $\br \in \mathfrak{T}_\delta$ one has that
$\mathsf{G}_{\br,\IC}:L^2(\D) \rightarrow L^2(\D)$ is linear and bounded.
For each $n\in \IN$ and $\br \in \mathfrak{T}_\delta$, set
\begin{equation}
	\mathsf{g}^{(n)}_{\br,\IC}
	(\bx,\by)
	\coloneqq
	n \nu
	\mathsf{p}_{\br,\IC}(\bx,\by)
	\chi'
	\left(
		n\norm{\bx-\by}^\nu
	\right)
	\norm{\bx-\by}^{\nu-2}
	(\bx-\by).
\end{equation}
Hence, for each $n\in \IN$ and $\br \in \mathfrak{T}_\delta$
and using item (ii) in the assumptions of \cref{eq:holomorphic_L2},
which is also part of the assumptions of the current result, we obtain
\begin{equation}
	\snorm{\mathsf{g}^{(n)}_{\br,\IC}}
	\leq
	n
	\nu
	C(\delta,\mathfrak{T})
	\snorm{
	\chi'
	\left(
		n\norm{\bx-\by}^\nu
	\right)
	}
	\norm{\bx-\by}^{-1},
	\quad
	(\bx,\by)
	\in
	\left(\D \times \D\right)^\star.
\end{equation}
Hence, 
\begin{equation}
\begin{aligned}
	\norm{
		\nabla_{\bx}
		{\mathsf{P}}^{(n)}_{\br,\IC}
		-
		\mathsf{G}_{\br,\IC}
	}_{\mathscr{L}\left(L^1(\D), L^1(\D)\right)}
	\leq
	&
	\underset{\by \in {\D}}{\operatorname{ess} \sup} 
	\int\limits_{\D}
	\snorm{
		1
		-
		\chi
		\left(
			n\norm{\bx-\by}^\nu
		\right)
	}
	\snorm{
		\nabla_{\bx}
		\mathsf{p}_{\br,\IC}(\bx,\by)
	}
	\dd \bx
	\\
	&
	+
	\underset{\by \in {\D}}{\operatorname{ess} \sup}
	\int\limits_{\D} 
	\norm{
		\mathsf{g}^{(n)}_{\br,\IC}
		(\bx,\by)
	}
	\dd \bx.
\end{aligned}
\end{equation}
Observe that since $\chi'(t) = 0$ for $t\geq 1$
one has
\begin{equation}\label{eq:der_P_conv_1}
\begin{aligned}
	\underset{\by \in {\D}}{\operatorname{ess} \sup}
	\int\limits_{\D} 
	\norm{
		\mathsf{g}^{(n)}_{\br,\IC}
	}
	\dd \bx
	&
	\leq
	n
	\nu
	C(\delta,\mathfrak{T})
	\,
	\underset{\by \in {\D}}{\operatorname{ess} \sup}
	\int\limits_{\D}
	\snorm{
	\chi'
	\left(
		n\norm{\bx-\by}^\nu
	\right)
	}
	\norm{\bx-\by}^{-1}
	\dd \bx
	\\
	&
	\leq
	n
	\nu
	C(\delta,\mathfrak{T})
	\,
	\underset{\by \in {\D}}{\operatorname{ess} \sup}
	\int\limits_{\D  \cap B(\by,n^{-\frac{1}{\nu}})}
	\frac{
		\mathrm{d}\bx
	}{
		\norm{\bx -\by}
	}
	\dd \bx
	\\
	&
	\leq
	n
	\nu
	C(\delta,\mathfrak{T})
	\omega_d
	\,
	\int\limits_{0}^{n^{-\frac{1}{\nu}}}
	\varrho^{d-2}
	\dd \varrho
	\\
	&
	\leq
	n
	\nu
	C(\delta,\mathfrak{T})
	\omega_d
	\,
	\left.
		\frac{
			\varrho^{d-1}
		}{
			d-1
		}
	\right\rvert_{0}^{n^{-\frac{1}{\nu}}}
	=
	\nu
	\frac{
		C(\delta,\mathfrak{T})
	}{
		d-1
	}
	n^{-\frac{d-1}{\nu}+1},
\end{aligned}
\end{equation}
where we have used that $\nu \in [0,d-1)$ and $C(\delta,\mathfrak{T})$ is as in \cref{eq:singular_kernel}.

Next, recalling \cref{eq:bound_kernel}, we estimate
{\small
\begin{equation}\label{eq:der_P_conv_2}
\begin{aligned}
	\underset{\by \in {\D}}{\operatorname{ess} \sup} 
	\int\limits_{\D}
	\snorm{
		1
		-
		\chi
		\left(
			n\norm{\bx-\by}^\nu
		\right)
	}
	\snorm{
		\nabla_{\bx}
		\mathsf{p}_{\br,\IC}(\bx,\by)
	}
	\dd \bx
	&
	\leq
	C(\delta,\mathfrak{T})
	\;
	\underset{\by \in {\D}}{\operatorname{ess} \sup} 
	\int\limits_{\D  \cap B(\by,n^{-\frac{1}{\nu}})}
	\frac{
		\mathrm{d}\bx
	}{
		\norm{\bx -\by}^{\zeta}
	}
	\dd \bx
	\\
	&
	\leq
	\frac{
		C(\delta,\mathfrak{T})
	}{
		d-\zeta
	}
	\omega_d
	n^{-\frac{d-\zeta}{\nu}}.
\end{aligned}
\end{equation}
}%
Consequently, it follows from \cref{eq:der_P_conv_1}
and \cref{eq:der_P_conv_2} that
\begin{equation}
	\sup_{\br \in \mathfrak{T}_\delta}
	\norm{
		\nabla_{\bx} {\mathsf{P}}^{(n)}_{\br,\IC}
		-
		\mathsf{G}_{\br,\IC}
	}_{\mathscr{L}\left(L^1(\D), L^1(\D)\right)}
	\leq
	\widetilde{C}(\delta,\mathfrak{T})
	\omega_d
	\left(
		\frac{
			n^{-\frac{d-1}{\nu}+1}
		}{
			d-1
		}
		+
		\frac{
			n^{\frac{d-\zeta}{\nu}}
		}{
			d-\zeta
		}
	\right).
\end{equation}
The same bound holds for 
$\norm{
		\nabla_{\bx} {\mathsf{P}}^{(n)}_{\br,\IC}
		-
		\mathsf{G}_{\br,\IC}
	}_{\mathscr{L}\left(L^\infty(\D), L^\infty(\D)\right)}$.
Then, recalling the Riesz-Thorin interpolation theorem stated in 
\Cref{sec:rt_theorem} we get
\begin{equation}
	\lim_{n \rightarrow \infty}
	\sup_{\br \in \mathfrak{T}_\delta}
	\norm{
		\nabla_{\bx}
		{\mathsf{P}}^{(n)}_{\br,\IC}	
		-
		\mathsf{G}_{\br,\IC}
	}_{\mathscr{L}\left(L^2(\D), L^2(\D)\right)}
	=
	0.
\end{equation}
Consequently, for each $\br \in \mathfrak{T}_\delta$
the operator ${\mathsf{P}}_{\br,\IC}:L^2(\D) \rightarrow H^1(\D)$
is linear and bounded, and for $\varphi \in L^2(\D)$
\begin{equation}
	\nabla_{\bx}
	\left(
		{\mathsf{P}}_{\br,\IC}
		\,
		\varphi
	\right)
	({\bf x})
	=
	\int\limits_{\D}
	\nabla_{\bx}
	\mathsf{p}_{\br,\IC}(\bx,\by)
	\varphi(\by)
	\dd{\by},
	\quad
	\text{a.e. }
	{\bx}
	\in 
	{\D}.	
\end{equation}
In the proof of \cref{eq:holomorphic_L2} we have shown that for each
$n\in \IN$ the map
$\mathfrak{T}_\delta \ni \br \mapsto \mathsf{p}^{(n)}_{\br,\IC} \in L^\infty(\D \times\D)$
is holomorphic.

To conclude the proof, it suffices to show that the map
\begin{equation}
	\mathfrak{T}_\delta
	\ni
	\br
	\mapsto
	\nabla_{\bx}
	{\mathsf{P}}_{\br,\IC}
	\in
	\mathscr{L}\left(L^2(\D), L^2(\D)\right).
\end{equation}
is holomorphic. To do so, we may argue as in the proof of \Cref{eq:holomorphic_L2},
thus concluding the proof of the current result. 

\end{proof}	

\subsection{Domain Transformations}
Let $\widehat\D$ be a simply connected, bounded Lipschitz domain
in $\IR^d$, $d\in \{2,3\}$.
We consider a subset $\mathfrak{T}$
of $W^{1,\infty}(\widehat\D,\IR^d)$, $d\in \{2,3\}$, of bijective and bi-Lipschitz
transformations such that for each 
$\br \in \mathfrak{T}$ the set $\D_\br= \br(\widehat\D)$ defines a 
simply connected, bounded Lipschitz domain in $\IR^d$.

\begin{assumption}\label{assump:transformations_T}
The set $\mathfrak{T}$ is a compact subset of
$W^{1,\infty}(\widehat\D,\IR^d)$, $d\in \{2,3\}$. 
For each $\br \in \mathfrak{T}$, 
$\br^{-1} \in W^{1,\infty}(\D_\br,\IR^d)$
and $\D_\br=\br(\widehat{\D})$ is a simply connected,
bounded Lipschitz domain.
\end{assumption}

\begin{assumption}\label{assump:piecewise_domains_T}
There exist $M \in \IN$ and 
$\widehat{\D}^{(1)},\dots,\widehat\D^{(M)}$ simply connected,
bounded Lipchitz domains in $\IR^d$, $d\in \{2,3\}$,
satisfying
\begin{equation}
	\widehat{\D}
	=
	\normalfont\text{int}
	\bigcup_{i=1}^{M}
	\overline{
		\widehat{\D}^{(i)}
	}
	\quad
	\text{and}
	\quad
	\widehat{\D}^{(i)} \cap \widehat{\D}^{(j)}
	= 
	\emptyset,
	\quad
	i,j \in \{1,\dots,M\},
	\quad
	i\neq j,
\end{equation}
such that for each $r\in \mathfrak{T}$ it holds
\begin{equation}
	\D_{\br}
	=
	\normalfont\text{int}
	\bigcup_{i=1}^{M}
	\overline{
		\D^{(i)}_{\br}
	}
	\quad
	\text{and}
	\quad
	\D^{(i)}_{\br} \cap \D^{(j)}_{\br}
	= 
	\emptyset,
	\quad
	i,j \in \{1,\dots,M\},
	\quad
	i\neq j,
\end{equation}
with $\D^{(i)}_{\br} \coloneqq \br\left(\widehat{\D}^{(i)}\right)$.
In addition, we assume that for each $\br\in \mathfrak{T}$
the set $\D^{(i)}_r$ defines a simply connected, bounded Lipschitz
domain in $\IR^d$, $d\in \{2,3\}$.
\end{assumption}

\begin{assumption}\label{assump:piecewise_analytic_kappa}
Consider the setting of \Cref{assump:piecewise_domains_T}.
We assume that there exist bounded Lipschitz domains $\D^{1}_{\normalfont\text{H}},
\dots,\D^{M}_{\normalfont\text{H}} \subset \IR^{d}$, $d\in \{2,3\}$,
satisfying
\[
	\bigcup_{\br \in \mathfrak{T}} 
	\D^{(i)}_\br
	\subset
	\overline{\D^{(i)}_{\text{H}}},
	\quad
	i\in \{1,\dots,M\}
\]
such that $\kappa$ admits an analytic extension into
an open neighborhood of $\D^{(i)}_{\text{H}}$. 
\end{assumption}

\begin{remark}
\Cref{assump:piecewise_domains_T,assump:piecewise_analytic_kappa}
have been introduced 
with the goal of allowing $\beta$ to be piecewise analytic.
This allows to include possibly discontinuous $\kappa$ defined
over a suitable partition of the working domain.
\end{remark}

For each $\br \in \mathfrak{T}$ and $u \in L^2(\widehat{\text{D}})$
we define 
$\tau_{\br}$ as $
	\left(
		\tau_{\br}u
	\right)(\widehat\bx)
	\coloneqq
	(u\circ \br)(\widehat\bx),
	\quad
	\bx
	\in 
	\widehat{\D}.
$
This operator may be extended to a scale of Sobolev spaces
Let $\mathfrak{T}$ be a set of admissible domain transformations.
Then for each $\br \in \mathfrak{T}$
\begin{align}\label{eq:pullback_operator}
	L^2(\D_\br)
	\ni
	u
	\mapsto
	\tau_{\br}u
	\in 
	L^2(\widehat\D)
	\quad
	\text{and}
	\quad
	H^1(\D_\br)
	\ni
	u
	\mapsto
	\tau_{\br}u
	\in 
	H^1(\widehat\D)
\end{align}
define isomorphisms, i.e.
$\tau_{\br} \in \mathscr{L}\left(L^2(\D_\br),L^2(\widehat\D)\right), 
\mathscr{L}\left(H^1(\D_\br),H^1(\widehat\D)\right)$ for each 
$\br \in \mathfrak{T}$.

For $\varepsilon>0$ and $d\in \{2,3\}$ let us define
\begin{equation}
	\mathfrak{T}_\varepsilon
	\coloneqq
	\left \{
		\br \in W^{1,\infty}(\widehat\D;\IC^d):
		\;
		\exists
		\,
		\widetilde{\br}
		\in 
		\mathfrak{T}
		\;\;
		\text{such that}
		\;\;
		\norm{\br-\widetilde{\br}}_{ W^{1,\infty}(\widehat{\D};\IC^d)}
		<
		\varepsilon
	\right \}.
\end{equation}	
and
$
	\left(\widehat{\D} \times \widehat{\D}\right)^\star 
	\coloneqq 
	\left\{
		(\widehat\bx,\widehat\by) 
		\in \widehat{\D} \times \widehat{\D}: 
		\widehat\bx\neq\widehat\by 
	\right\}.
$

The following technical results asserts the existence
of an open \emph{complex} neighborhood of a set 
of admissible domain transformations $\mathfrak{T}$
within which certain properties are preserved.

\begin{proposition}\label{eq:prop_existence_delta}
Let $\mathfrak{T}$ be a set of domain transformations
satisfying \Cref{assump:transformations_T}.
There exists $\delta = \delta(\mathfrak{T})>0$
and a constant $\zeta(\delta,\mathfrak{T})>0$, depending on $\delta$
and $\mathfrak{T}$ only, such that
\[\label{eq:}
	\inf_{\br \in {\mathfrak{T}_\delta}}
	\inf_{(\widehat\bx,\widehat\by) \in (\widehat\D \times \widehat\D)^\star }
	\Re
	\left \{
		\frac{(\br(\widehat\bx)-\br(\widehat\by))\cdot ((\br(\widehat\bx)-\br(\widehat\by)))}{\norm{\widehat\bx-\widehat\by}^2}
	\right\}
	\geq
	\zeta(\delta,\mathfrak{T})
	>0.
\]
\end{proposition}

\begin{proof}
Exactly as in the proof of \cite[Proposition 3.2]{henriquez2021shape}.
\end{proof}

As a consequence of \cite[Theorem 3.8]{CSZ18}
we have the holomorphic dependence of certain 
relevant quantities upon domain transformations.

\begin{lemma}[{\cite[Theorem 3.8]{CSZ18}}]\label{lmm:holomorphy_jacobian}
Let $\mathfrak{T}$ be a set of domain transformations
satisfying \Cref{assump:transformations_T}.
Then there exist $\varepsilon>0$ such that the map\footnote{
Observe that $D$ denotes the Jacobian matrix whereas $\widehat{\text{D}}$
refers to the reference domain.}
\begin{equation}
	\mathcal{J}:
	\mathfrak{T}
	\rightarrow
	L^{\infty}(\widehat\D):
	\br
	\mapsto
	J_{\br}
	\coloneqq
	\det D \br
\end{equation}
admits a bounded holomorphic extension
into $\mathfrak{T}_\varepsilon$.
\end{lemma}

\subsection{Shape Holomorphy of the Lippmann-Schwinger VIE}
\label{sec:shape_hol_volumen_operator}
Throughout, let $\mathfrak{T}$ be a set of domain transformations
satisfying \Cref{assump:transformations_T}. 
As in \Cref{sec:vol_int_eq}, we consider $\kappa_0 \in \IR_+$
and, provided that \cref{assump:transformations_T,assump:piecewise_analytic_kappa}
are satisfied, we consider
$\kappa \in L^{\infty}(\D_\text{H})$, where 
\begin{equation}
	\D_\text{H} 
	\coloneqq
	\bigcup_{i=1}^{M}\D^{(i)}_\text{H}.
\end{equation}

For each $\br \in \mathfrak{T}$ and
$u \in H^{1}(\D_{\br})$, we define 
\begin{equation}
	\left(
		\mathsf{A}_{\br} u
	\right)
	(\bx)
	\coloneqq
	\int\limits_{\D_\br}
	G^{(\kappa_0)}(\bx-\by)
	\beta(\by)
	u(\by)
	\dd\by,
	\quad
	\bx\in \D_\br,
\end{equation}
which corresponds to the operator
$\mathsf{A}$ introduced in \cref{eq:VIO}
on the domain $\D_\br$, i.e. we explicitely
introduce the dependence on the domain transformation
$\br \in \mathfrak{T}$ through a subscript.
For each $\br \in \mathfrak{T}$, we set 
\begin{equation}\label{eq:operator_A_r}
	\widehat{\mathsf{A}}_{\br}
	\coloneqq
	\tau_{\br}
	\;
	\mathsf{A}_{\br}
	\;
	\tau^{-1}_{\br}
	\in 
	\mathscr{L}
	\left(	
		X
		,
		X
	\right),
\end{equation}
where the $X$ in \cref{eq:operator_A_r}
can be either $H^1(\widehat{\D})$ or $L^2(\widehat{\D})$,
and the stated mapping properties follow
straightforwardly from \cref{eq:VIO} and
\cref{eq:pullback_operator}.

The operator $\widehat{\mathsf{A}}_{\br}$ 
defined in \cref{eq:operator_A_r} admits
the following explicit expression: For
$\widehat{u} \in H^1(\widehat{\D})$ 
and
$\br \in \mathfrak{T}$ 
one has
\begin{equation}
	\left(
		\widehat{\mathsf{A}}_{\br}
		\widehat{u}
	\right)
	(\widehat\bx)
	\coloneqq
	\int\limits_{\widehat{\D}}
	\mathsf{g}_{\br}^{(\kappa_0)}
	\left(
		\widehat\bx,\widehat\by
	\right)
	\beta(\br(\widehat\by))
	\widehat{u}(\widehat\by)
	\dd
	\widehat\by,
	\quad
	\widehat\bx\in\widehat{\D},
\end{equation}
where
\begin{equation}\label{eq:transformerd_kernel}
	\mathsf{g}_{\br}^{(\kappa_0)}
	\left(
		\widehat\bx,\widehat\by
	\right)
	\coloneqq
	G^{(\kappa_0)}(\br(\widehat\bx)-\br(\widehat\by))
	J_{\br}(\widehat\by),
	\quad
	(\widehat\bx,\widehat\by)
	\in
	(\widehat{\D} \times \widehat{\D})^\star,
\end{equation}
with $G^{(\kappa_0)}$ as in \cref{eq:green_G_omega}.

We are interested in the smoothness 
of the \emph{domain-to-operator} map
\begin{equation}
	\mathfrak{T}
	\ni
	\br 
	\mapsto 
	\widehat{\mathsf{A}}_{\br}
	\in
	\mathscr{L}
	\left(	
		X
		,
		X
	\right),
\end{equation}
where, again, $X$
denotes either $H^1(\widehat{\D})$ or $L^2(\widehat{\D})$.
We start by considering the latter case.

\begin{theorem}\label{thm:shape_hol_VIO_L2}
Let $\mathfrak{T}$ be a set of domain transformations
satisfying \Cref{assump:transformations_T},
and let \Cref{assump:piecewise_domains_T,assump:piecewise_analytic_kappa}
be satisfied.
Then there exists $\delta>0$ such that the map
\begin{equation}\label{eq:map_hol_A}
	\mathcal{A}:
	\mathfrak{T}
	\rightarrow
	\mathscr{L}
	\left(	
		L^2(\widehat{\D})
		,
		L^2(\widehat{\D})
	\right):
	\br
	\mapsto
	\widehat{\mathsf{A}}_{\br}
\end{equation}
admits a bounded holomorphic extension into
$\mathfrak{T}_\delta$.
\end{theorem}

\begin{proof}
As we aim to allow piecewise analytic coefficients,
we consider the following decomposition:
Under \Cref{assump:piecewise_domains_T}, 
one may write for each $\br \in \mathfrak{T}$ 
and $\widehat{u} \in L^2(\widehat{\D})$
\[
	\left(
		\widehat{\mathsf{A}}_{\br}
		\,
		\widehat{u}
	\right)
	(\widehat{{\bf x}})
	=
	\sum_{i=1}^{M}
	\int\limits_{\widehat{\D}^{(i)}}
	\mathsf{g}_{\br}^{(\kappa_0)}
	\left(
		\widehat{\bx},\widehat{\by}
	\right)
	\widehat{\beta}_{\br}(\widehat{\by})
	\widehat{u}(\widehat{\by})
	\dd
	\widehat{\by},
	\quad
	\widehat{\bx}
	\in 
	\widehat{\D},
\]
where $\widehat{\beta}_{\br} \coloneqq \beta\circ \br$.
In addition, we set for $\widehat{u} \in L^2(\widehat{\D})$ and $i\in \{1,\dots,M\}$
\begin{align}
	\left(
		\widehat{\mathsf{A}}^{(i)}_{\br}
		\,
		\widehat{u}
	\right)
	(\widehat{\bx})
	=
	\int\limits_{\widehat{\D}^{(i)}}
	\mathsf{g}_{\br}^{(\kappa_0)}
	\left(
		\widehat{\bx},\widehat{\by}
	\right)
	\widehat{\beta}_{\br}(\widehat{\by})
	\widehat{u}(\widehat{\by})
	\dd
	\widehat{\by},
	\quad
	\widehat{\bx}
	\in 
	\widehat{\D},
\end{align}
thus yielding
\begin{align}\label{eq:decomp_A}
	\widehat{\mathsf{A}}_{\br}
	=
	\sum_{i=1}^{M}
	\widehat{\mathsf{A}}^{(i)}_{\br}.
\end{align}
For each $i\in \{1,\dots,M\}$ and $\br \in \mathfrak{T}$ 
the operator $\widehat{\mathsf{A}}^{(i)}_{\br}: L^2(\widehat{\D})\rightarrow
L^2(\widehat{\D})$ is linear and bounded.
It follows from \Cref{eq:prop_existence_delta}
and \Cref{assump:piecewise_analytic_kappa} that there exists
$\delta>0$ such that for each $\br \in \mathfrak{T}_\delta$
and $\widehat{u}\in L^2(\widehat{\D})$ one can define
the following extension to complex-valued domain transformations:
\begin{equation}
	\left(
		\widehat{\mathsf{A}}^{(i)}_{\br,\IC}
		\,
		\widehat{u}
	\right)
	(\widehat{{\bf x}})
	\coloneqq
	\int\limits_{\widehat{\D}^{(i)}}
	\mathsf{g}_{\br,\IC}^{(\kappa_0)}
	\left(
		\widehat{\bx},\widehat{\by}
	\right)
	\widehat{\beta}_{\br}(\widehat{\by})
	\widehat{u}(\widehat{\by})
	\dd
	\widehat{\by},
	\quad
	\widehat{\bx}
	\in 
	\widehat{\D},
\end{equation}
where
\begin{equation}\label{eq:def_g_kappa_0}
	\mathsf{g}_{\br,\IC}^{(\kappa_0)}
	\left(
		\widehat{\bx},\widehat{\by}
	\right)
	=
	\left\{
	\begin{array}{cl}
	\frac{\imath}{4}
	\text{H}^{(1)}_0
	\left(
		\kappa_0
		\norm{\br(\widehat{\bx})-\br(\widehat{\by})}_\IC
	\right)
	J_{\br,\IC}(\widehat{\by}),  & d=2 \\
	\displaystyle
	\frac{
		\exp
		\left(
			\imath\kappa_0 
			\norm{\br(\widehat{\bx})-\br(\widehat{\by})}_\IC
		\right)	
	}{
		4\pi \norm{\br(\widehat{\bx})-\br(\widehat{\by})}_\IC
	}
	J_{\br,\IC}(\widehat{\by}), & d=3
	\end{array}
	\right.
	,
	\quad
	(\widehat{\bx},\widehat{\by}) 
	\in 
	\left(
		\widehat{\D}\times \widehat{\D}
	\right)^\star.
\end{equation}
Above, $\norm{{\boldsymbol{v}}}_\IC 
\coloneqq \sqrt{\boldsymbol{v}\cdot \boldsymbol{v}}$
for any $\boldsymbol{v} \in \mathfrak{U}$, $d\in \IN$,
and $\sqrt{\cdot}: \mathbb{C}\backslash (-\infty,0]\rightarrow \IC$
signifies the principal branch of the square root.

We verify item-by-item the assumptions of
\cref{eq:holomorphic_L2} for the case $d=2$.
\begin{itemize}
	\item[(i)]
	According to \cref{eq:prop_existence_delta} 
	there exists $\delta>0$ such that
	for each $\br \in \mathfrak{T}_\delta$ and
	$(\widehat\bx,\widehat\by) \in (\widehat\D \times \widehat\D)^\star$ one 
	has
	\[\label{eq:lowe_bound}	
		\Re
		\left\{
			(\br(\widehat\bx)-\br(\widehat\by))\cdot (\br(\widehat\bx)-\br(\widehat\by))
		\right\}
		\geq
		\zeta(\delta,\mathfrak{T})
		\norm{\widehat\bx-\widehat\by}^2
		>0.
	\]
	The Hankel function $\text{H}^{(1)}_0$ of the first kind of order zero is 
	analytic in in right complex half-plane, i.e. for $z\in \IC$ with $\Re{\{z\}}>0$.
	Furthermore, it follows from \cref{eq:prop_existence_delta} that
	for each $\br \in \mathfrak{T}_\delta$
	\[\label{eq:lower_bound_Q}
		\Re
		\left\{
			\frac{
			\norm{\br(\widehat{\bx})-\br(\widehat{\by})}_\IC	
			}{
			\norm{\widehat\bx-\widehat\by}
			}
		\right\}
		\geq
		\sqrt{
		\zeta(\delta,\mathfrak{T})
		}
		>0,
	\]
	Hence, recalling also \cref{lmm:holomorphy_jacobian},
	for $\delta >0$ as in \cref{eq:prop_existence_delta},
	each $(\widehat\bx,\widehat\by) \in (\widehat\D \times \widehat\D)^\star$ 
	that for $d=2$ the map
	\begin{equation}
		\mathfrak{T}_\delta
		\ni
		\br
		\mapsto
		\mathsf{g}_{\br,\IC}^{(\kappa_0)}
		\left(
			\widehat{\bx},\widehat{\by}
		\right)
		\in
		\IC
	\end{equation}
	is holomorphic.
	\item[(ii)]
	The function $\text{H}^{(1)}_0$ admits for $z\in \IC$ with $\Re{\{z\}}>0$ 
	the following representation (see e.g. \cite[Eq. 12.31, pp. 330]{Ober12})
	\begin{equation}
		\text{H}^{(1)}_0(z)
		=
		-\frac{2\imath}{\pi}
		\int\limits_{0}^{\infty}
		\frac{\exp{((\imath-t) z)}}{\sqrt{t^2-2\imath t}}
		\dd t,
	\end{equation}
where the branch cut of $\sqrt{t^2-2\imath t}$
is selected so that $\Re \{\sqrt{t^2-2\imath t} \}\geq 0$
for $t>0$.
Also, for $z\in \IC$ with $\Re{\{z\}}>0$ it holds
\begin{equation}\label{eq:upper_bound_Hankel}
\begin{aligned}
	\snorm{\text{H}^{(1)}_0(z)}
	\leq
	\frac{2}{\pi}
	\snorm{\exp({\imath z})}
	\int\limits_{0}^{\infty}
	\frac{\snorm{\exp{(-tz)}}}{\snorm{\sqrt{t^2-2\imath t}}}
	\dd
	t
	&
	\leq
	\frac{2}{\pi}
	\snorm{\exp({\imath z})}
	\int\limits_{0}^{\infty}
	\frac{\exp{(-t\Re{\{z\}})}}{\sqrt{2t}}
	\dd
	t \\
	&
	=
	\frac{\sqrt{2}}{\pi}
	\snorm{\exp({\imath z})}
	\int\limits_{0}^{\infty}
	\frac{\exp{(-t\Re{\{z\}})}}{\sqrt{t}}
	\dd
	t. 
\end{aligned}
\end{equation}
It follows that
\[
	\snorm{\text{H}^{(1)}_0(z)}
	\leq
	\snorm{\exp({\imath z})}
	\sqrt{\frac{2}{\pi}}
	\left.
	\frac{\text{erf}(\sqrt{\Re{\{z\}}} \sqrt{t})}{\sqrt{\Re{\{z\}}}}
	\right|_{0}^{\infty} 
	=
	\snorm{\exp({\imath z})}
	\sqrt{\frac{2}{\pi{\Re{\{z\}}}}},
\]
where 
$
	\text{erf}(t)
	=
	\frac{2}{\sqrt{\pi}}
	\int\limits_{0}^{t}
	\exp(-s^2)
	\dd
	s
$
corresponds to the \emph{error function}.
Observe that for each $\br \in \mathfrak{T}_\delta$
and $(\widehat{\bx},\widehat{\by}) 
	\in 
	\left(
		\widehat{\D}\times\widehat{\D}
	\right)^\star$
\begin{equation}
	\snorm{
	\mathsf{g}^{(\kappa_0)}_{\br,\IC}
	\left(
		\widehat{\bx},\widehat{\by}
	\right)
	}
	\leq
	\frac{1}{4\sqrt{\norm{\widehat\bx-\widehat\by}}}
	\sqrt{\frac{2}{\pi}}
	\frac{
	\snorm{
	\exp
	\left(
		\kappa_0 
		\norm{\br(\widehat{\bx})-\br(\widehat{\by})}_\IC
	\right)
	}
	}{
	\sqrt{
	\Re
	\left\{
		\frac{
		\norm{\br(\widehat{\bx})-\br(\widehat{\by})}_\IC	
		}{
		\norm{\widehat\bx-\widehat\by}
		}
	\right\}
	}
	}.
\end{equation}
For each $\br \in \mathfrak{T}_\delta$ there exists 
$\widetilde{\br} \in \mathfrak{T}$ such that
$\norm{\br -\widetilde{\br} }_{W^{1,\infty}(\widehat\D,\IR^d)}<\delta$, thus
\begin{equation}\label{eq:map_Bounded}
\begin{aligned}
	\snorm{
	\exp
	\left(
		\kappa_0 
		\norm{\br(\widehat{\bx})-\br(\widehat{\by})}_\IC
	\right)
	}
	&
	\leq
	\exp
	\left(
		\kappa_0 
		\text{diam}(\widehat\D)
		\norm{\br-\widetilde{\br}}_{W^{1,\infty}(\widehat\D,\IR^d)}
		+
		\norm{\widetilde{\br}}_{W^{1,\infty}(\widehat\D,\IR^d)}
	\right) 
	\\
	&
	\leq
	\exp
	\left(
		\kappa_0 
		\text{diam}(\widehat\D)
		\delta
	\right)
	\sup_{\widetilde{\br} \in \mathfrak{T}}
	\exp
	\left(
		\kappa_0 
		\text{diam}(\widehat\D)
		\norm{\widetilde{\br}}_{W^{1,\infty}(\widehat\D,\IR^d)}
	\right)
	\\
	&
	<
	\infty,
\end{aligned}
\end{equation}
where the boundedness of the last term follows from the continuity of
the map 
\begin{equation}
	\mathfrak{T} 
	\ni 
	\br 
	\mapsto 	
	\exp
	\left(
		\kappa_0 
		\text{diam}(\widehat\D)
		\norm{{\br}}_{W^{1,\infty}(\widehat\D,\IR^d)}
	\right)
	\in
	\IR
\end{equation}
and the compactness of $\mathfrak{T}$. 
Finally, recalling \cref{eq:lower_bound_Q} we obtain
that item (ii) in \cref{eq:holomorphic_L2} holds, i.e.
\[
	\snorm{
	\mathsf{g}^{(\kappa_0)}_{\br,\IC}
	\left(
		\widehat{\bx},\widehat{\by}
	\right)
	}
	\leq
	\frac{
		C(\delta,\mathfrak{T})
	}{
		\norm{\widehat\bx-\widehat\by}^{\half}
	},
	\quad
	(\widehat\bx,\widehat\by) \in (\widehat{\D} \times \widehat{\D})^\star.
\]
for some constant $C(\delta,\mathfrak{T})>0$ that depends
on $\delta$ and $\mathfrak{T}$ only.
\end{itemize}

Next, we verify item-by-item the assumptions of
\cref{eq:holomorphic_L2} for the case $d=3$.
\begin{itemize}
	\item[(i)]
	As in the two dimensional case, for each 
	$(\widehat\bx,\widehat\by) \in (\widehat\D \times \widehat\D)^\star$
	the map
	\begin{equation}
		\mathfrak{T}_\delta
		\ni
		\br
		\mapsto
		\norm{\br(\widehat{\bx})-\br(\widehat{\by})}_\IC
		\in
		\IC
	\end{equation}
	is holomorphic. Then, for 
	$d=3$ and recalling \cref{eq:lowe_bound},
	for each $(\widehat\bx,\widehat\by) \in (\widehat\D \times \widehat\D)^\star$
	the map
	\begin{equation}
		\mathfrak{T}_\delta
		\ni
		\br
		\mapsto
		\mathsf{g}_{\br,\IC}^{(\kappa_0)}
		\left(
			\widehat{\bx},\widehat{\by}
		\right)
	\end{equation}
	is holomorphic as well.
	\item[(ii)]
	Observe that for each $\mathfrak{T}_\delta$
	and each $(\widehat\bx,\widehat\by) \in (\widehat\D \times \widehat\D)^\star$
	\begin{align}
		\snorm{
		\mathsf{g}^{(\kappa_0)}_{\br,\IC}
		\left(
			\widehat{\bx},\widehat{\by}
		\right)
		}
		\leq
		\frac{1}{4\pi\norm{\widehat\bx-\widehat\by}}
		\frac{
	\snorm{
	\exp
	\left(
		\kappa_0 
		\norm{\br(\widehat{\bx})-\br(\widehat{\by})}_\IC
	\right)
	}
	}{
	\Re
	\left\{
		\frac{
		\norm{\br(\widehat{\bx})-\br(\widehat{\by})}_\IC	
		}{
		\norm{\widehat\bx-\widehat\by}
		}
	\right\}
	},
	\quad
	(\widehat{\bx},\widehat{\by}) 
	\in 
	(
		\widehat{\D}\times\widehat{\D}
	)^\star.
\end{align}
Recalling \cref{eq:lower_bound_Q} and \cref{eq:map_Bounded}
we may conclude that there exists $C(\kappa_0,\delta,\mathfrak{T})>0$,
depending on $\kappa_0$, $\delta$, and $\mathfrak{T}$, such that
\begin{equation}
	\snorm{
		\mathsf{g}^{(\kappa_0)}_{\br,\IC}
		\left(
		\widehat{\bx},\widehat{\by}
		\right)
	}
	\leq
	\frac{
		C(\kappa_0,\delta,\mathfrak{T})
	}{
		\norm{\widehat\bx-\widehat\by}
	},
	\quad
	(\widehat{\bx},\widehat{\by}) 
	\in 
	(
		\widehat{\D}\times\widehat{\D}
	)^\star.
\end{equation}
Thus, for the three dimensional case item (ii) in
\cref{eq:holomorphic_L2} holds as well.
\end{itemize}

Consequently, for either $d=2$ or $d=3$, the map in \cref{eq:map_hol_A}
admits a bounded holomorphic extension to $\mathfrak{T}_\delta$ for some $\delta>0$.
\end{proof}

%

Finally, we present the corresponding shape holomorphy result when considering
the VIOs as elements of $\mathscr{L}\left(	L^2(\widehat{\D}),H^1(\widehat{\D})\right)$.

\begin{theorem}
Let $\mathfrak{T}$ be a set of domain transformations
satisfying \Cref{assump:transformations_T},
and let \Cref{assump:piecewise_domains_T,assump:piecewise_analytic_kappa}
be satisfied.
Then there exists $\delta>0$ such that the map
\begin{align}\label{eq:map_hol_A_H1}
	\mathcal{A}:
	\mathfrak{T}
	\rightarrow
	\mathscr{L}
	\left(	
		L^2(\widehat{\D})
		,
		H^1(\widehat{\D})
	\right):
	\br
	\mapsto
	\widehat{\mathsf{A}}_{\br}
\end{align}
admits a bounded holomorphic extension into
$\mathfrak{T}_\delta$.
\end{theorem}

\begin{proof}
As in the proof of \Cref{thm:shape_hol_VIO_L2} we consider the
decomposition stated in \cref{eq:decomp_A}.
We proceed to verify the assumptions of \Cref{eq:holomorphic_H1}.

We verify item-by-item the assumptions of
\cref{eq:holomorphic_H1} for the case $d=2$ item-by-item.
\begin{itemize}
\item[(i)]  Arguing as in the proof of \Cref{thm:shape_hol_VIO_L2}, in particular item (i) for $d=2$,
we may conclude that there exists $\delta>0$ such that 
for each $(\widehat{\bx},\widehat{\by}) \in \left(\widehat{\D}\times \widehat{\D}\right)^\star$ the map 
\begin{equation}
	\mathfrak{T}_\delta \ni \br
	\mapsto
	\nabla_{\widehat\bx}
	\;
	\mathsf{g}_{\br,\IC}^{(\kappa_0)}
	\left(
		\widehat{\bx},\widehat{\by}
	\right)
\end{equation}
is holomorphic.
\item[(ii)] For $(\widehat{\bx},\widehat{\by}) \in \left(\widehat{\D}\times \widehat{\D}\right)^\star$, we have that
\begin{equation}
	\nabla_{\widehat\bx}
	\;
	\mathsf{g}_{\br,\IC}^{(\kappa_0)}
	\left(
		\widehat{\bx},\widehat{\by}
	\right)
	=
	-
	\kappa_0
	\frac{\imath}{4}
	\frac{
		\left(
			\br(\widehat{\bx})-\br(\widehat{\by})
		\right)^{\top}
		\cdot
		\left(D\boldsymbol{r}\right)(\widehat{\bx})
	}{
		\norm{\br(\widehat{\bx})-\br(\widehat{\by})}_\IC
	}
	\text{H}^{(1)}_1
	\left(
		\kappa_0
		\norm{\br(\widehat{\bx})-\br(\widehat{\by})}_\IC
	\right)
	J_{\br,\IC}(\widehat{\by})
\end{equation}
The function $\text{H}^{(1)}_1$ admits for $z\in \IC$ with $\Re{\{z\}}>0$ 
the following representation (see e.g. \cite[Eq. 12.31, pp. 330]{Ober12})
\begin{equation}
		\text{H}^{(1)}_1(z)
		=
		\frac{2\imath}{\pi}
		\int\limits_{0}^{\infty}
		\frac{(\imath-t)\exp{((\imath-t) z)}}{\sqrt{t^2-2\imath t}}
		\dd t,
\end{equation}
As in step (ii) of the proof of \cref{thm:shape_hol_VIO_L2} for the two dimensional case,
we may show that for any $z\in \IC$ with $\Re{\{z\}}>0$ it holds
\begin{equation}
	\snorm{
		\text{H}^{(1)}_1(z)
	}
	\leq
	\snorm{\exp({\imath z})}
	\left(
		\sqrt{\frac{2}{\pi{\Re{\{z\}}}}}
		+
		\frac{2}{\pi{\Re{\{z\}}}}
	\right).
\end{equation}
Consequently, for $(\widehat{\bx},\widehat{\by}) 
	\in 
	\left(
		\widehat{\D}\times \widehat{\D}
	\right)^\star$
\begin{equation}
	\snorm{
	\nabla_{\widehat\bx}
	\;
	\mathsf{g}_{\br,\IC}^{(\kappa_0)}
	\left(
		\widehat{\bx},\widehat{\by}
	\right)
	}
	\leq
	\frac{
	C(\mathfrak{T},\delta,\kappa_0)
	}{
	\norm{\widehat\bx-\widehat\by}
	}
	\left(
	\frac{
		\norm{\widehat\bx-\widehat\by}^{\frac{1}{2}}
	}{
	\sqrt{
	\Re
	\left\{
		\frac{
		\norm{\br(\widehat{\bx})-\br(\widehat{\by})}_\IC	
		}{
		\norm{\widehat\bx-\widehat\by}
		}
	\right\}
	}
	}
	+
	\frac{
	1
	}{
	{
	\Re
	\left\{
		\frac{
		\norm{\br(\widehat{\bx})-\br(\widehat{\by})}_\IC	
		}{
		\norm{\widehat\bx-\widehat\by}
		}
	\right\}
	}
	}
	\right),
\end{equation}
where $C(\mathfrak{T},\delta,\kappa_0)$ depends on $\mathfrak{T}$, $\delta$, and $\kappa_0$.
\end{itemize}

Now, we verify item-by-item the assumptions of
\cref{eq:holomorphic_H1} for the case $d=3$ item-by-item.
\begin{itemize}
\item[(i)]  Arguing as in the proof of \Cref{thm:shape_hol_VIO_L2}, in particular item (i) for $d=3$,
we may conclude that there exists $\delta>0$ such that 
for each $(\widehat{\bx},\widehat{\by}) \in \left(\widehat{\D}\times \widehat{\D}\right)^\star$ the map 
\begin{equation}
	\mathfrak{T}_\delta \ni \br
	\mapsto
	\nabla_{\widehat\bx}
	\;
	\mathsf{g}_{\br,\IC}^{(\kappa_0)}
	\left(
		\widehat{\bx},\widehat{\by}
	\right)
\end{equation}
is holomorphic.
\item[(ii)] For $(\widehat{\bx},\widehat{\by}) \in \left(\widehat{\D}\times \widehat{\D}\right)^\star$ we have
\begin{equation}
\begin{aligned}
	\nabla_{\widehat\bx}
	\mathsf{g}_{\br,\IC}^{(\kappa_0)}
	\left(
		\widehat{\bx},\widehat{\by}
	\right)
	=
	&
	\frac{\imath\kappa_0}{4\pi}
	\frac{
		\left(
			\br(\widehat{\bx})-\br(\widehat{\by})
		\right)^{\top}
		\cdot
		\left(D\boldsymbol{r}\right)(\widehat{\bx})
	}{
		\norm{\br(\widehat{\bx})-\br(\widehat{\by})}_\IC
	}
	\frac{
		\exp
		\left(
			\imath\kappa_0 
			\norm{\br(\widehat{\bx})-\br(\widehat{\by})}_\IC
		\right)	
	}{
		\norm{\br(\widehat{\bx})-\br(\widehat{\by})}_\IC
	}
	J_{\br,\IC}(\widehat{\by})
	\\
	&
	-
	\frac{	
		1
	}{
		4\pi
	}
	\exp
	\left(
		\imath\kappa_0 
		\norm{\br(\widehat{\bx})-\br(\widehat{\by})}_\IC
	\right)	
	\frac{
		\left(
			\br(\widehat{\bx})-\br(\widehat{\by})
		\right)^{\top}
		\cdot
		\left(D\boldsymbol{r}\right)(\widehat{\bx})
	}{
		\norm{\br(\widehat{\bx})-\br(\widehat{\by})}^3_\IC
	}
	J_{\br,\IC}(\widehat{\by}).
\end{aligned}
\end{equation}
Thus, for $(\widehat{\bx},\widehat{\by}) \in \left(\widehat{\D}\times \widehat{\D} \right)^\star$
\begin{equation}
	\snorm{
		\nabla_{\widehat\bx}
		\mathsf{g}_{\br,\IC}^{(\kappa_0)}
		\left(
			\widehat{\bx},\widehat{\by}
		\right)
	}
	\leq
	\frac{
	C(\mathfrak{T},\delta,\kappa_0)
	}{
	\norm{\widehat\bx-\widehat\by}^2
	}
	\left(
	\frac{
		\norm{\widehat\bx-\widehat\by}
	}{
	\Re
	\left\{
		\frac{
		\norm{\br(\widehat{\bx})-\br(\widehat{\by})}_\IC	
		}{
		\norm{\widehat\bx-\widehat\by}
		}
	\right\}
	}
	+
	\frac{
	1
	}{
	{
	\Re
	\left\{
		\frac{
		\norm{\br(\widehat{\bx})-\br(\widehat{\by})}_\IC	
		}{
		\norm{\widehat\bx-\widehat\by}
		}
	\right\}^2
	}
	}
	\right),
\end{equation}
where, again, $C(\mathfrak{T},\delta,\kappa_0)$ depends on $\mathfrak{T}$, $\delta$, and $\kappa_0$.
\end{itemize}

Consequently, for either $d=3$ or $d=3$ the map in \cref{eq:map_hol_A_H1}
admits a bounded holomorphic extension to $\mathfrak{T}_\delta$ for some $\delta>0$.

\end{proof}

\subsection{Shape Holomorphy of the Domain-to-Solution map}
In \Cref{sec:shape_hol_volumen_operator}, based on the results proved in \Cref{sec:holomorphic_volume_operators}, we have established the holomorphic dependence 
of the VIO appearing in the Lippmann-Schwinger VIE.
\begin{corollary}\label{cor:shape_hol_d2s_map}
Let $\mathfrak{T}$ be a set of domain transformations
satisfying \Cref{assump:transformations_T},
and let \Cref{assump:piecewise_domains_T,assump:piecewise_analytic_kappa} be satisfied.
In addition, assume that $\uinc$ is entire in $\mathbb{R}^d$, $d \in \{2,3\}$.
Then there exists $\delta>0$ such that the map
\begin{align}
	\mathcal{S}:
	\mathfrak{T}
	\rightarrow
	X:
	\br
	\mapsto
	\widehat{u}_{\br}
	\coloneqq
	\left(\mathsf{Id} - \widehat{\mathsf{A}}_{\br}\right)^{-1} \widehat{u}^{\normalfont\text{inc}}_{{\br}} ,
\end{align}
where $X \in \{L^2(\widehat{\normalfont\text{D}}),H^1(\widehat{\normalfont\text{D}})\}$ 
and $\widehat{u}^{\normalfont\text{inc}}_{{\br}} \coloneqq \tau_{\br} \uinc$
admits a bounded holomorphic extension into
$\mathfrak{T}_\delta$.
\end{corollary}

\begin{proof}
This results follows from the well-posedness of \Cref{LSvar} in the physical domain $\D_{\br}$
for each $\br\in \mathfrak{T}$, the fact that $\tau_{\br} \in \mathscr{L}_{\text{iso}}\left(H^1(\D_\br),H^1(\widehat\D)\right)$ for each $\br \in \mathfrak{T}$, the properties of $\uinc$, and \cite[Proposition 4.20]{HS21}.
\end{proof}

However, in practical applications we do not have direct access to the
domain-to-solution map but to a suitable approximation in a finite dimensional
subspace. To this end, we consider the following Galerkin formulation in a finite dimensional
subspace $V_h \subset H^1(\widehat{\normalfont\text{D}})$ as described in \Cref{sec:num_dict}: 
For each $\br \in \mathfrak{T}$ we seek $\widehat{u}_{h,\br} \in V_h$ such that
\begin{equation}\label{eq:bilinearLS_hatD}
	\widehat{\mathsf{a}}_\br(\widehat{u}_{h,\br}, \widehat{v}_h) 
	\coloneqq 
	\dual{\widehat{u}_{h,\br}}{\widehat{v}_h}_{\widehat{\text{D}}} 
	- 
	\dual{\widehat{\mathsf{A}}_{\br} \widehat{u}_{h,\br}}{\widehat{v}_h}_{\widehat{\text{D}}}
	= 
	\dual{\normalfont\uinc}{\widehat{v}_h}_{\widehat{\text{D}}},
\end{equation}
holds for all $\widehat{v}_h\in V_h$.

The variational problem stated in \cref{eq:bilinearLS_hatD} is well-posed. 
Indeed, one can readily notice that
$\widehat{\mathsf{A}}_{\br} \in  \mathscr{L}\left(L^2(\widehat{\normalfont\text{D}}),H^1(\widehat{\normalfont\text{D}})\right)$ is compact, and together with \Cref{th:discrete-inf-sup} and Fredholm's alternative 
one may conclude that for each $\br \in \mathfrak{T}$ there exists $h_0(\br)$ \emph{depending on $\br \in \mathfrak{T}$} such that for any $0<h<h_0(\br)$ the aforementioned variational problem is well-posed.

The next result address the holomorphic dependence of $\widehat{u}_{h,\br} \in V_h$
upon the domain transformation $\br \in\mathfrak{T}$.

\begin{corollary}\label{cor:shape_hol_d2s_map_discrete}
Let $\mathfrak{T}$ be a set of domain transformations
satisfying \Cref{assump:transformations_T},
and let \Cref{assump:piecewise_domains_T,assump:piecewise_analytic_kappa} be satisfied.
In addition, assume that $\uinc$ is entire in $\mathbb{R}^d$, $d \in \{2,3\}$.
Then, there exists $\delta>0$ and $h_0= h_0(\mathfrak{T},\delta)$ (depending only on
$\mathfrak{T}$ and $\delta$) such that for $0<h<h_0$ the map
\begin{align}
	\mathcal{S}_h:
	\mathfrak{T}
	\rightarrow
	V_h:
	\br
	\mapsto
	\widehat{u}_{h,\br},
\end{align}
admits a bounded holomorphic extension into
$\mathfrak{T}_\delta$.
\end{corollary}

\begin{proof}
The proof follows the same steps as in \Cref{cor:shape_hol_d2s_map}. 
However, it must be noticed that $h_0$ and $\delta>0$ in this result, unlike to that of the well-posednes
argument presented before, do not depend on the particular instance of the domain transformation $\br \in \mathfrak{T}$. This can be proved by using a finite covering argument and the compactness of $\mathfrak{T}$
as a subset of $W^{1,\infty}(\widehat{\text{D}};\mathbb{R}^d)$, $d\in \{2,3\}$, cf. \Cref{assump:transformations_T}.
\end{proof}

\section{Applications to Forward and Inverse UQ}
\label{sec:applications}
In this section, we describe the implications of the results presented in \Cref{sec:vol_int_eq,sec:shape_holomorphy} into Forward and Inverse UQ.
\subsection{Parametric Holomorphy}
In concrete applications, domain deformations are parametrically
defined through a sequence $\y \in \mathbb{U} \coloneqq [-1,1]^\IN$, i.e.
we could in principle afford an infinite, yet countable, number of parameters. 

For $\varrho>1$, we consider the Bernstein ellipse in the complex plane
\begin{align}
	\mathcal{E}_s
	\coloneqq 
	\left\{
		\frac{z+z^{-1}}{2}: \; 1\leq \snorm{z}\leq s
	\right \} 
	\subset \IC.
\end{align}
This ellipse has foci at $z=\pm 1$ and semi-axes of length 
$a\coloneqq  (s+s^{-1})/2$ and $b \coloneqq  (s-s^{-1})/2$.
Let us consider the tensorized poly-ellipse
\begin{align}
\mathcal{E}_{\boldsymbol{\rho}} \coloneqq  \bigotimes_{j\geq1} \mathcal{E}_{\rho_j} \subset \IC^{\mathbb{N}},
\end{align}
where $\boldsymbol\rho \coloneqq  \{\rho_j\}_{j\geq1}$ is such that $\rho_j>1$, for $j\in \mathbb{N}$.

\begin{definition}[{\cite[Definition 2.1]{CCS15}}]\label{def:bpe_holomorphy}
Let $X$ be a complex Banach space equipped with the norm $\norm{\cdot}_{X}$. 
For $\varepsilon>0$ and $p\in(0,1)$, we say that map 
$\mathbb{U}  \ni  \y \mapsto  u(\y)  \in  X$
is $(\boldsymbol{b},p,\varepsilon)$-holomorphic if and only if:
\begin{enumerate}
	\item\label{def:bpe_holomorphy1}
	The map $\mathbb{U} \ni {\y} \mapsto u(\y) \in X$ is uniformly bounded.
	\item\label{def:bpe_holomorphy2}
	There exists a positive sequence $\boldsymbol{b}\coloneqq (b_j)_{j\geq 1} \in \ell^p(\mathbb{N})$ 
	and a constant $C_\varepsilon>0$ such that for any sequence 
	$\boldsymbol\rho\coloneqq (\rho_j)_{j\geq1}$ 
	of numbers strictly larger than one that is
	$(\boldsymbol{b},\varepsilon)$-admissible, i.e.~satisfying
	$\sum_{j\geq 1}(\rho_j-1) b_j  \leq  \varepsilon$,
	the map $\y \mapsto u(\y)$ admits a complex
	extension $\boldsymbol{z} \mapsto u(\boldsymbol{z})$ 
	that is holomorphic with respect to each
	variable $z_j$ on a set of the form 
	\begin{align}
		\mathcal{O}_{\boldsymbol\rho} 
		\coloneqq  
		\displaystyle{\bigotimes_{j\geq 1}} \, \mathcal{O}_{\rho_j},
	\end{align}
	where $ \mathcal{O}_{\rho_j}=
	\{z\in\IC\colon\operatorname{dist}(z,[-1,1])<\rho_j-1\}$.
	\item
	This extension is bounded on $\mathcal{E}_{\boldsymbol\rho}$ according to
	\begin{equation}
		\sup_{\boldsymbol{z} \in \mathcal{E}_{\boldsymbol{\rho}}} 
		\norm{u(\boldsymbol{z})}_{X}  \leq C_\varepsilon.
	\end{equation}
\end{enumerate}
\end{definition}

\subsection{Affine-Parametric Domain Transformations}
\label{sec:affine_parametric_transformations}
In the two-dimensional case, we consider the affine-parametric domain
transformations introduced in \cite[Section 5.3]{CSZ18}.
Let $\widehat{\varrho} \in W_{\text {per }}^{1, \infty}(0,2 \pi)$
and  set
\begin{equation}\label{eq:radius_y}
	 \varrho_{\boldsymbol{y}}(\varphi)
	 \coloneqq \widehat{\varrho}(\varphi)+\sum_{j \in \mathbb{N}} y_{j}
	 \varrho_{j}(\varphi) \in W_{\mathrm{per}}^{1, \infty}(0,2 \pi).
\end{equation}
Let us set $b_j = \left\|\varrho_{j}\right\|_{W_{\text {per}}^{1, \infty}(0,2 \pi)}$
and throughout we assume that
$\boldsymbol{b} = \{b_j\}_{j \in \mathbb{N}} \in \ell^p(\mathbb{N}) \subset \ell^1(\mathbb{N}) $
for some $0<p<1$ together with
\begin{equation}
	\varrho_{\min } 
	\leq 
	\widehat{\varrho}(\varphi)-\sum_{j \in \mathbb{N}}\left|\varrho_{j}(\varphi)\right| 
	\leq 
	\widehat{\varrho}(\varphi)+\sum_{j \in \mathbb{N}}\left|\varrho_{j}(\varphi)\right| 
	\leq 
	\varrho_{\max }
\end{equation}
for some $0<\varrho_{\min } < \varrho_{\max }<\infty$.
Set
\begin{equation}\label{eq:domains_param_ref}
\begin{aligned}
	&\widehat{\mathrm{D}}
	\coloneqq
	\left\{\vartheta(\cos \varphi, \sin \varphi)^{\top}: 0 \leq \vartheta \leq \widehat{\varrho}(\varphi),
	\quad \varphi \in [0,2\pi)\right\}, \\
	&\mathrm{D}_\y
	\coloneqq
	\left\{\vartheta(\cos \varphi, \sin \varphi)^{\top}: 0 \leq \vartheta \leq \varrho_\y(\varphi),
	\quad \varphi \in [0,2\pi)\right\},
	\quad
	\y \in \mathbb{U},
\end{aligned}
\end{equation}
and for each $\y \in \mathbb{U}$ consider $\br_\y: \widehat{\text{D}} \rightarrow \text{D}_\y: \widehat{\bx} \mapsto {\frac{\varrho_\y(\varphi)}{\widehat{\varrho}(\varphi)}\widehat{\bx}} $.
Setting
\begin{equation}
	\bar{\br}(\widehat{\bx} )\coloneqq
	\widehat{\bx} 
	\quad \text {and} 
	\quad 
	{\br}_{j}(\widehat{\bx})
	\coloneqq
	\frac{\varrho_{j}(\varphi)}{\widehat{\varrho}(\varphi)}\widehat{\bx},
\end{equation}
where $\widehat{\bx} = \sigma (\cos \varphi,\sin \varphi)^\top$, we have
\begin{equation}\label{eq:transform_r_y}
	\br_\y(\widehat{\bx})
	=
	\bar{\br}(\widehat{\bx})
	+
	\sum_{j\geq1}
	y_j
	{\br}_{j}(\widehat{\bx}),
	\quad
	\widehat{\bx}\in \widehat{\text{D}},
	\quad
	\y = \{y_j\}_{j\in \mathbb{N}} 
	\in \mathbb{U}.
\end{equation}
The set 
\begin{equation}\label{eq:set_T}
	\mathfrak{T} \coloneqq \left \{\br_\y: \; \y \in \mathbb{U} \right \}
\end{equation}
is a compact subset of $W^{1,\infty}(\widehat{\text{D}};\mathbb{R}^{d})$,
$d \in \{2,3\}$, and as it has been established in 
\cite[Lemma 5.7]{CSZ18}. 

Equipped with this, we can establish the following result.
\begin{corollary}\label{cor:bps_hol_sol_map}
Suppose that $\mathfrak{T}$ as in \cref{eq:set_T} fullfils  
\Cref{assump:transformations_T,assump:piecewise_domains_T,assump:piecewise_analytic_kappa}.
In addition, assume that $\boldsymbol{b} \in \ell^p(\mathbb{N})$ for some $0<p<1$.
The following maps are $(\boldsymbol{b},p,\varepsilon)$-holomorphic for some $\epsilon>0$
(independent of $h$ in item (ii)):
\begin{itemize}
	\item[(i)]
	{\sf Parameter-to-Solution map}
	$
		\mathcal{S}:
		\mathbb{U}
		\rightarrow
		X:
		\y
		\mapsto
		\widehat{u}_{\y}
		\coloneqq 
		\widehat{u}_{\br_\y}.
	$
	\item[(ii)]
	{\sf Discrete Parameter-to-Solution map:}
	For any $0<h<h_0$,
	$
		\mathcal{S}_h:
		\mathbb{U}
		\rightarrow
		X:
		\y
		\mapsto
		\widehat{u}_{h,\y} \coloneqq \widehat{u}_{h,\br_\y}.
	$
\end{itemize}

\end{corollary}

\begin{proof}
This result follows straightforwardly from \Cref{cor:shape_hol_d2s_map,cor:shape_hol_d2s_map_discrete}
together with \cite[Lemma 5.8]{CSZ18}.
\end{proof}

\subsection{Forward UQ: Expected Value of the Point evaluation}
\label{sec:forward_UQ}
Let $\Omega$ be a bounded Lipschitz domain with boundary $\partial \Omega$.
For each $\y \in \mathbb{U}$, let $\widehat{u}_{h,\y}$ be as in \Cref{cor:bps_hol_sol_map}.
Using the representation formula stated in \cref{eq:LS-rep}
\begin{equation}\label{eq:QoI}
	{u}_{h,\y}({\bf x})
	=
	\uinc ({\bf x})
	+ 
	\N^{(\kappa_0)}_{\y}(\beta {u}_{h,\y})({\bf x}),
	\quad
	{\bf x} 
	\in \partial \Omega,
\end{equation}
where ${u}_{h,\y} = \tau^{-1}_{\y}\widehat{u}_{h,\y}$. 
We are interested in the computation of the expected value of 
the QoI defined in \cref{eq:QoI} with respect to the uniform measure in the 
parameter space.

Observe that
\begin{equation}
\begin{aligned}
	\left(
		\N^{(\kappa_0)}_\y(\beta {u}_{h,\y})
	\right)
	(\bx)
	&
	=
	\int\limits_{\D_\y}
	G^{(\kappa_0)}(\bx-\by)
	\beta(\by)
	{u}_{h,\y}(\by)
	\dd\by
	\\
	&
	=
	\int\limits_{\widehat{\D}}
	\mathsf{g}_{\y}^{(\kappa_0)}
	\left(
		\bx,\widehat\by
	\right)
	\beta(\br_\y(\widehat\by))
	\widehat{u}_{h,\y}(\widehat\by)
	\dd
	\widehat\by,
	\quad
	\bx\in\partial \Omega,
\end{aligned}
\end{equation}
where, for each $\y \in \mathbb{U}$, $\mathsf{g}_{\y}^{(\kappa_0)}$ is defined as
\begin{equation}
	\mathsf{g}_{\y}^{(\kappa_0)}
	\left(
		{\bx},\widehat{\by}
	\right)
	=
	\left\{
	\begin{array}{cl}
	\frac{\imath}{4}
	\text{H}^{(1)}_0
	\left(
		\kappa_0
		\norm{\bx-\br_\y(\widehat{\by})}
	\right)
	J_{\y}(\widehat{\by}),  & d=2 \\
	\displaystyle
	\frac{
		\exp
		\left(
			\imath\kappa_0 
			\norm{\bx-\br_\y(\widehat{\by})}
		\right)	
	}{
		4\pi \norm{\bx-\br_\y(\widehat{\by})}
	}
	J_{\y}(\widehat{\by}), & d=3
	\end{array}
	\right.
	,
	\quad
	({\bx}, \widehat{\by}) 
	\in 
	\partial \Omega \times \widehat{\D},
\end{equation}
where $J_{\y} \coloneqq J_{\br_\y}$, $\y \in \mathbb{U}$.

\begin{corollary}\label{cor:bpe_hol_QoI}
Suppose that $\mathfrak{T}$ as in \cref{eq:set_T} fullfils  
\Cref{assump:transformations_T,assump:piecewise_domains_T,assump:piecewise_analytic_kappa}.
In addition, assume that $\boldsymbol{b} \in \ell^p(\mathbb{N})$ for some $0<p<1$.
For each $ \bx \in \partial \Omega$, the map 
\begin{equation}\label{eq:map_QoI_hol}
	\mathcal{F}:
	\mathbb{U} \rightarrow \IC: \y \mapsto 
	\left(
		\N^{(\kappa_0)}_\y(\beta {u}_{h,\y})
	\right)
	(\bx)
\end{equation}
is $(\boldsymbol{b},p,\varepsilon)$-holomorphic for some $\epsilon>0$.
\end{corollary}

\begin{proof}
As in the proof of \Cref{thm:shape_hol_VIO_L2} since for each $\bx \in \partial \Omega$
one has 
\begin{equation}
	\inf_{\y \in \mathbb{U}} \inf_{\hat{\bx} \in \widehat{\D}} \norm{\bx - \br_\y(\hat{\bx})} >0
\end{equation}
it holds that the map 
$
	\mathbb{U}
	\ni
	\y
	\mapsto
	\norm{\bx - \br_\y(\hat{\bx})} 
	\in 
	\mathbb{R}
$
is $(\boldsymbol{b},p,\varepsilon)$-holomorphic for some $\varepsilon>0$.

Next, it follows from \Cref{assump:transformations_T,assump:piecewise_domains_T,assump:piecewise_analytic_kappa}
together with \Cref{lmm:holomorphy_jacobian} that for each $\bx \in \partial \Omega$ the map 
$\mathbb{U} \ni \y \mapsto \mathsf{g}_{\y}^{(\kappa_0)}
	\left(
		\bx,\cdot
	\right)
	\beta(\br_\y(\cdot)) \in L^\infty(\widehat{\D})$
	is $(\boldsymbol{b},p,\varepsilon)$-holomorphic as well. It follows
	from \Cref{cor:bps_hol_sol_map} that
	$\y\mapsto \widehat{u}_{h,\y}$ is so as well.
	Therefore, we may conclude that $\mathcal{Q}$ as in \cref{eq:map_QoI_hol}
	also has this property, thus concluding the proof.
\end{proof}

For a given $\bx \in \partial \Omega$, we are interested 
in computing the mean value of $\mathcal{Q}$ over the parameter space $\mathbb{U}$,
i.e. 
\begin{equation}\label{eq:integral_QoI}
	\mathcal{I}
	\coloneqq
	\int\limits_{\mathbb{U}}
	\mathcal{F}(\y)
	\mu(\dd \y),
\end{equation}
where $\mu = \bigotimes_{j\in \mathbb{N}} \frac{\lambda}{2}$ is the product measure
with $\frac{\lambda}{2}$ being $\frac{1}{2}$ the Lebesgue measure in $[-1,1]$.
To this end, we firstly consider only the first $s \in \mathbb{N}$ dimensions in the parameter
space, and consider instead
\begin{equation}\label{eq:integral_QoI_Approx}
	\mathcal{I}_s
	\coloneqq
	\int\limits_{\mathbb{U}^{(s)}}
	\mathcal{F}(\y_{\{1:s\}})
	\mu^{(s)}(	\dd \y). 
\end{equation}
as an approximation of \cref{eq:integral_QoI}, with
$\mu^{(s)} = \displaystyle\bigotimes_{j=1}^{s} \frac{\lambda}{2}$ and $\mathbb{U}^{(s)} = [-1,1]^s$,
and for $\y \in \mathbb{U}^{(s)}$ we set $\y_{\{1:s\}} = (y_1,y_2,\dots,y_s,0,0,0,\dots) \in \mathbb{U}$.

The effect of truncation in the parameter space when considering
the approximation \cref{eq:integral_QoI} by \cref{eq:integral_QoI_Approx}
produces an error that decays as $s^{-\left(\frac{1}{p}-1\right)}$, see e.g. 
\cite{DLC16}.
Let us set consider a $s$-dimensional, $N$ points, equal weights quadrature rule for the 
approximation of $\mathcal{I}_s$ defined in \cref{eq:integral_QoI_Approx}
\begin{equation}
	\mathcal{Q}_{N,s}
	\coloneqq
	\frac{1}{N}
	\sum_{i=1}^{N}
	\mathcal{F}\left(\y^{(i)}_{\{1:s\}}\right),
\end{equation}
where $\left\{\y^{(1)},\dots \y^{(N)}\right\}$ are the corresponding quadrature points
in $\mathbb{U}^{(s)}$.

As a consequence of \cref{cor:bpe_hol_QoI}, the 
numerical approximation of \cref{eq:integral_QoI_Approx} can 
be efficiently performed using HoQMC as in, e.g., \cite{DLC16}
with a provably rate of convergence $N^{-\frac{1}{p}}$, i.e.
\begin{equation}
	\snorm{
		\mathcal{I}_s
		-
		\mathcal{Q}_{N,s}
	}
	\lesssim
	N^{-\frac{1}{p}},
\end{equation}
with $0<p<1$ as in \cref{cor:bpe_hol_QoI}, and $N$ being
the total number of quadrature points, which in practice equates to 
the total number of forward model evaluations.
Notably, this rate of convergence is independent of the chosen 
dimension of truncation $s \in \mathbb{N}$ in the parameter space.

\subsection{Inverse UQ: Bayesian Shape Identification}
\label{eq:bayesian_shape_iden}
We rapidly recall the most import aspects of Bayesian inverse problems following
\cite{stuart2010inverse,SS12}.
Assume given a \emph{parameter--to--solution map} $G: \mathbb{U} \rightarrow X $
with $X$ being a separable Banach space over $\IC$ and a
\emph{prior} probability measure $\mu_0$ on $\mathbb{U}$. 
Our goal is to infrm the prior measure from observational data. In the
following we will assume that $\mu_0$ is the uniform probability
measure. In addition, we assume given an \emph{observation operator}
$O: X \rightarrow \mathbb{R}^K$. In turn, this defines the
\emph{uncertainty--to--observation} map as
$
	\mathcal{G}
	\coloneqq 
	O \circ G: 
	\mathbb{U} \rightarrow \mathbb{R}^K.
$
We assume observations $\Upsilon \in \IR^K$ following the following model
\begin{align}\label{eq:model_noise}
	\Upsilon = \mathcal{G}({\y^\star})+\eta,
\end{align} 
where $\eta \sim \mathcal{N}(0,\Sigma)$, $\Sigma \in \mathbb{R}^{K \times K}_\text{sym}$ is a symmetric, positive definite covariance matrix and $\y^\star$ is the \emph{ground truth}.

The following result addresses how the so-called posterior \emph{posterior probability measure} $\mu^{\Upsilon}$ can be expressed in terms of the prior one when the data $\Upsilon \in \IR^K$ is incorporated.

\begin{theorem}[{\cite[Theorem 2.1]{SS12}}]\label{thm:bayes_thm}
Assume that $\mathcal{J}: \mathbb{U} \rightarrow \IR^K$ is bounded and continuous.
Then $\mu^\Upsilon(\dd \y)$, the distribution of $\y \in \mathbb{U}$
given the data $\Upsilon \in \mathbb{R}^K$, is absolutely continuous with respect to
$\mu_0(\dd \y)$, i.e. there exists a parametric density $\Theta(\y)$
such that for each $\y \in \mathbb{U}$ the \emph{Radon-Nikodym} derivative
is given by
\begin{align}
	\frac{\dd \mu^\Upsilon}{\dd \mu_0}(\y) 
	= 
	\frac{1}{Z} 
	\Theta(\y)
\end{align}
with the posterior density
\begin{align}\label{eq:post_density}
	\Theta(\y) 
	\coloneqq \exp \left(-\Phi_{\Sigma}(\Upsilon, \boldsymbol{y}) \right)
	\quad
	\text{and}
	\quad
	Z
	\coloneqq 
	\int\limits_{\mathbb{U}}\Theta(\y) \mu_0(\dd \boldsymbol{y})>0,
\end{align}
where
\begin{align}
	\Phi_{\Sigma}(\Upsilon, \boldsymbol{y}) 
	= 
	\frac{1}{2}
	\left(
		\Upsilon - \mathcal{G}({\boldsymbol{y}}) 
	\right)^\top 
	\Sigma^{-1}  
	\left(
		\Upsilon - \mathcal{G}({\boldsymbol{y}})
	\right).
\end{align}
\end{theorem}

The applicability of this result is the context of Bayesian shape identification
is as follows. 
\begin{itemize}
	\item[(i)]
	The parameter-to-solution map is given by $G: \y \mapsto \widehat{u}_{h,\y}(\widehat\by)$.
	\item[(ii)]
	Consider a set of $K$ points $\{\bx_k\}_{k=1}^{K}$ in $\partial \Omega$. 
	The parameter-to-observation map then reads
	\begin{equation}\label{eq:param_to_measurement_map}
		\mathcal{G}:
		\y
		\mapsto 
		\begin{pmatrix}
			\Re\left\{ {u}_{h,\y}({\bf x}_1) \right\} \\
			\Im\left\{ {u}_{h,\y}({\bf x}_1) \right\} \\
			\vdots \\
			\Re\left\{ {u}_{h,\y}({\bf x}_K) \right\} \\
			\Im\left\{ {u}_{h,\y}({\bf x}_K) \right\} \\
		\end{pmatrix}
		\in 
		\mathbb{R}^{2K},
	\end{equation}
	where ${u}_{h,\y}({\bf x}_i)$ is as in \cref{eq:QoI}.
	\item[(iii)]
	We are interested in computing the expected value of the boundary, i.e. we consider the
	QoI $\phi: \y \mapsto \br_\y$ provided that the model is informed by
	some data $\Upsilon$.
	To this end, we recall \cref{thm:bayes_thm} and compute the expected value with 
	respect to the posterior distribution $\mu^\Upsilon(\dd \y)$ as follows
	\begin{align}\label{eq:expected_qoi}
		\mathbb{E}^{\mu^\Upsilon}[\phi]
		=
   		\int\limits_{\mathbb{U}}\phi(\y)\mu^\Upsilon(\dd\y)
    		=
		\frac{1}{Z}
		\int\limits_{\mathbb{U}}\phi(\y)\Theta(\y) \mu_0(\dd\boldsymbol{y}),
	\end{align}
	with $Z$ as in \cref{eq:post_density}.
\end{itemize}
As it has been discussed in \Cref{sec:forward_UQ}, the numerical 
computation of the high-dimensional integral appearing in \cref{eq:expected_qoi}
can be done using high-dimensional quadrature rules such as HoQMC.

\section{Numerical Results}
\label{sec:numerical_experiments}
We consider the reference domain $\widehat{\D}$ to be the unit disk centered at the origin.
The wavenumber $\kappa$ in \cref{eq:Helmholtz} is given by
\begin{equation}
	\kappa(\bx) \coloneqq \left\{\begin{array}{cl}  1, & \bx\in \IR^2 \setminus \overline{\widehat{\D}}, \\
								  2 + \tfrac{1}{2}\norm{\bx}^2, & \bx\in \widehat{\D}. 
								  \end{array}\right.
\end{equation}
We consider the setting described in \Cref{sec:affine_parametric_transformations}, i.e.
we study \Cref{LSvar} in a domain $\D_{\y} \coloneqq \br_\y(\widehat{\D})$ as in \cref{eq:domains_param_ref}
with $\br_\y:\widehat{\text{D}} \rightarrow \D_{\y} $ as in \cref{eq:transform_r_y}.
The radius is affine-parametric and given by a Karhunen-Lo\`eve-type expansion for each
$ \y\in \mathbb{U}$ of the form
\begin{equation}\label{eq:radius}
	\varrho_\y(\varphi) 
	\coloneqq 
	1 
	+ 
	\theta\sum_{j = 1}^{\infty} 
	j^{-\zeta} 
	\left( y_{2j} \cos(j\varphi ) + y_{2j-1} \sin(j\varphi) \right), 
	\quad \varphi\in [0, 2\pi],
\end{equation}
with $\zeta>1$, which is a concrete instance of \cref{eq:radius_y}.
Comparing with \cref{eq:radius_y} we get that
\begin{equation}
	\varrho_{2j}
	\coloneqq
	\theta
	j^{-\zeta}
	\cos(j\varphi ) 
	\quad
	\text{and}
	\quad
	\varrho_{2j-1}
	\coloneqq
	\theta
	j^{-\zeta}
	\sin(j\varphi ),
	\quad
	j \in \mathbb{N},
\end{equation}
and thus $\boldsymbol{b}= \{b_j\}_{j \in \mathbb{N}}$ with 
$b_j = \left\|\varrho_{j}\right\|_{W_{\text {per}}^{1, \infty}(0,2 \pi)}$ 
belongs to $\ell^p(\mathbb{N})$ for $p<\frac{1}{-1+\zeta}$.
Therefore, as discussed in \Cref{sec:forward_UQ}, the best possible rate of convergence is
 in the approximation of the expected value of a given QoI is $N^{-\frac{1}{p}}$ for any $p<\frac{1}{-1+\zeta}$.

We choose the parameters $\zeta \in \{2, 3\}$ and $\theta \in  \{\tfrac{1}{4}, \tfrac{3}{4}\}$ for our experiments. We truncate the sum in \cref{eq:radius} at $s = 100$. We study the convergence of QMC approximations to high-dimensional integrals relevant in the context of forward and inverse uncertainty quantification in this setting. We generate QMC points based on two different approaches: (a) randomly-shifted lattice rules (RLR); (b) interlaced-polynomial lattice rules (IPL) with interlacing parameter $\alpha$ (for further details we refer to
\cite{DLC16}).
Both can be generated by means of the \textsc{QMC4PDE}\footnote{\url{https://people.cs.kuleuven.be/~dirk.nuyens/qmc4pde/}} library. Relevant parameters are the asymptotic behaviour of the sequence $\boldsymbol{b} $. In all of our examples, reference solutions are computed with $N = 2048$ QMC points.

\subsection{Forward UQ: Expected Value of the Point evaluation}\label{sec:forward-uq}
\begin{figure}[ht!]
\centering
	\begin{subfigure}[b]{0.4\textwidth}
		\includegraphics[width=1.0\linewidth]{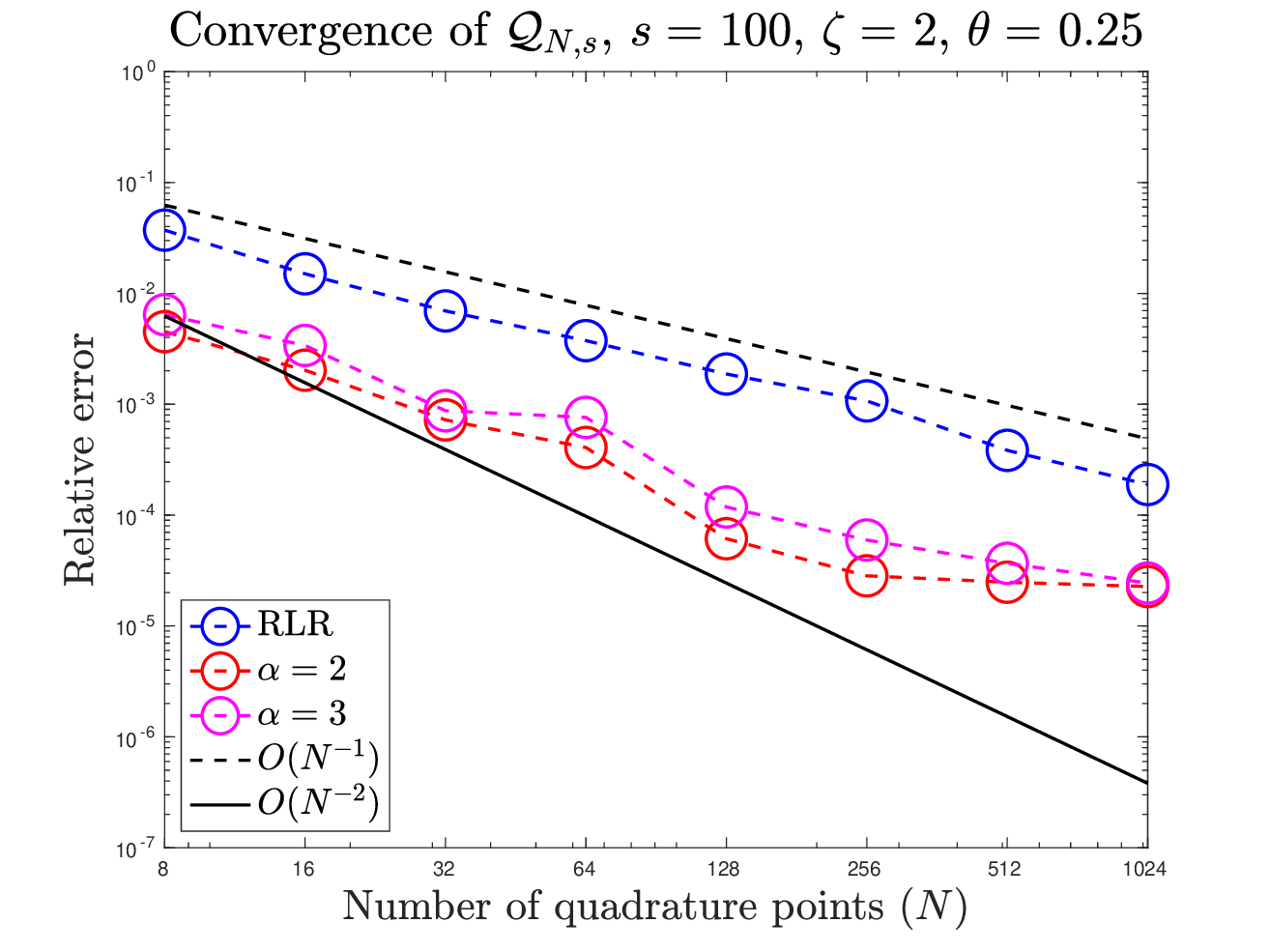}	
		\caption{$\zeta=2, \ s = 100$.}
	\end{subfigure}
	\begin{subfigure}[b]{0.4\textwidth}
		\includegraphics[width=1.0\linewidth]{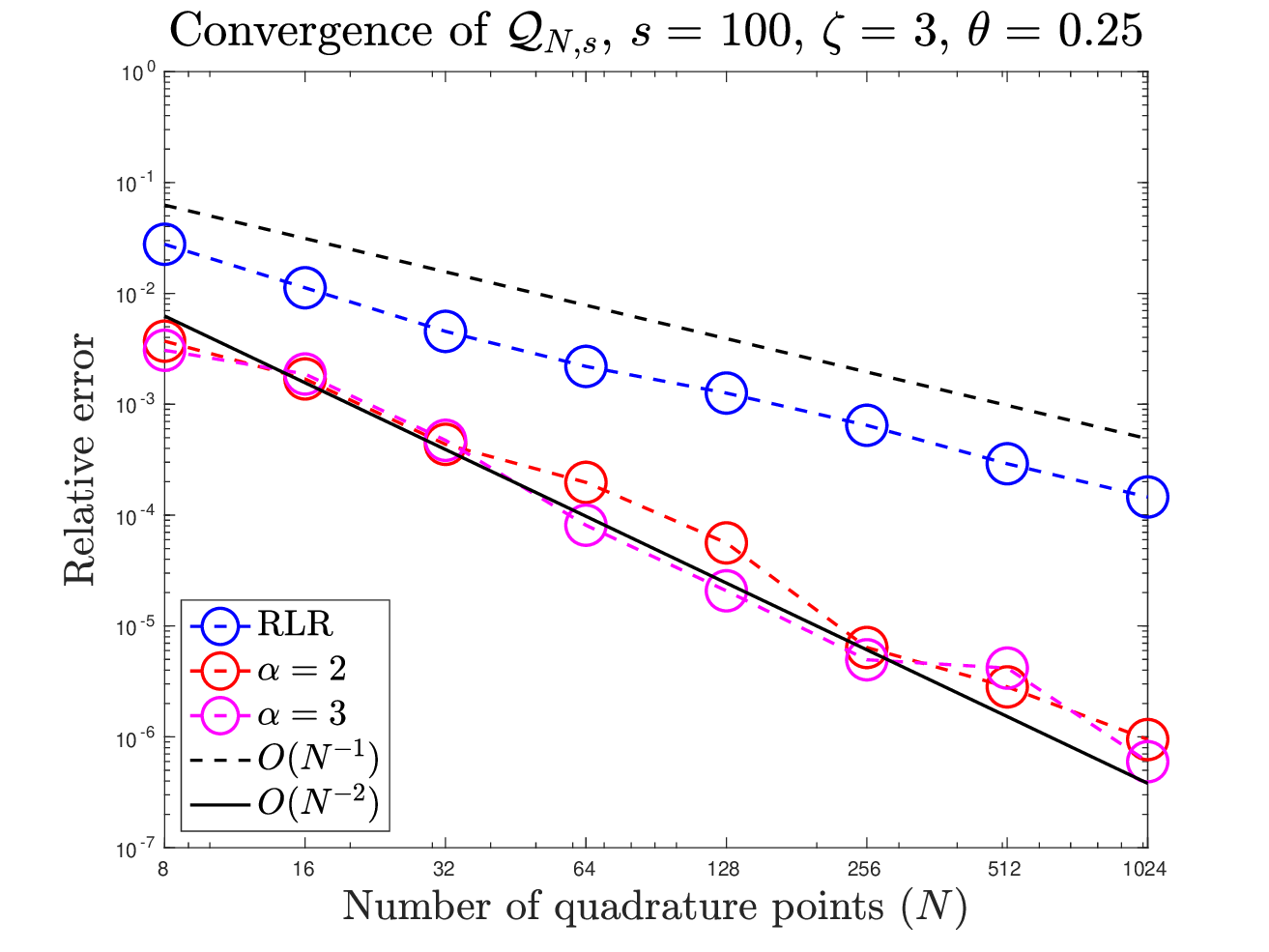}	
		\caption{$\zeta=3, \ s = 100$.}
	\end{subfigure}
	\caption{Convergence of $\mathcal{Q}_{N, s}$ vs. number of QMC points ($N$): $\theta = \tfrac{1}{4}$.}
	\label{fig:Test1-Forward}
\end{figure}

\begin{figure}[ht!]
\centering
	\begin{subfigure}[b]{0.4\textwidth}
		\includegraphics[width=1.0\linewidth]{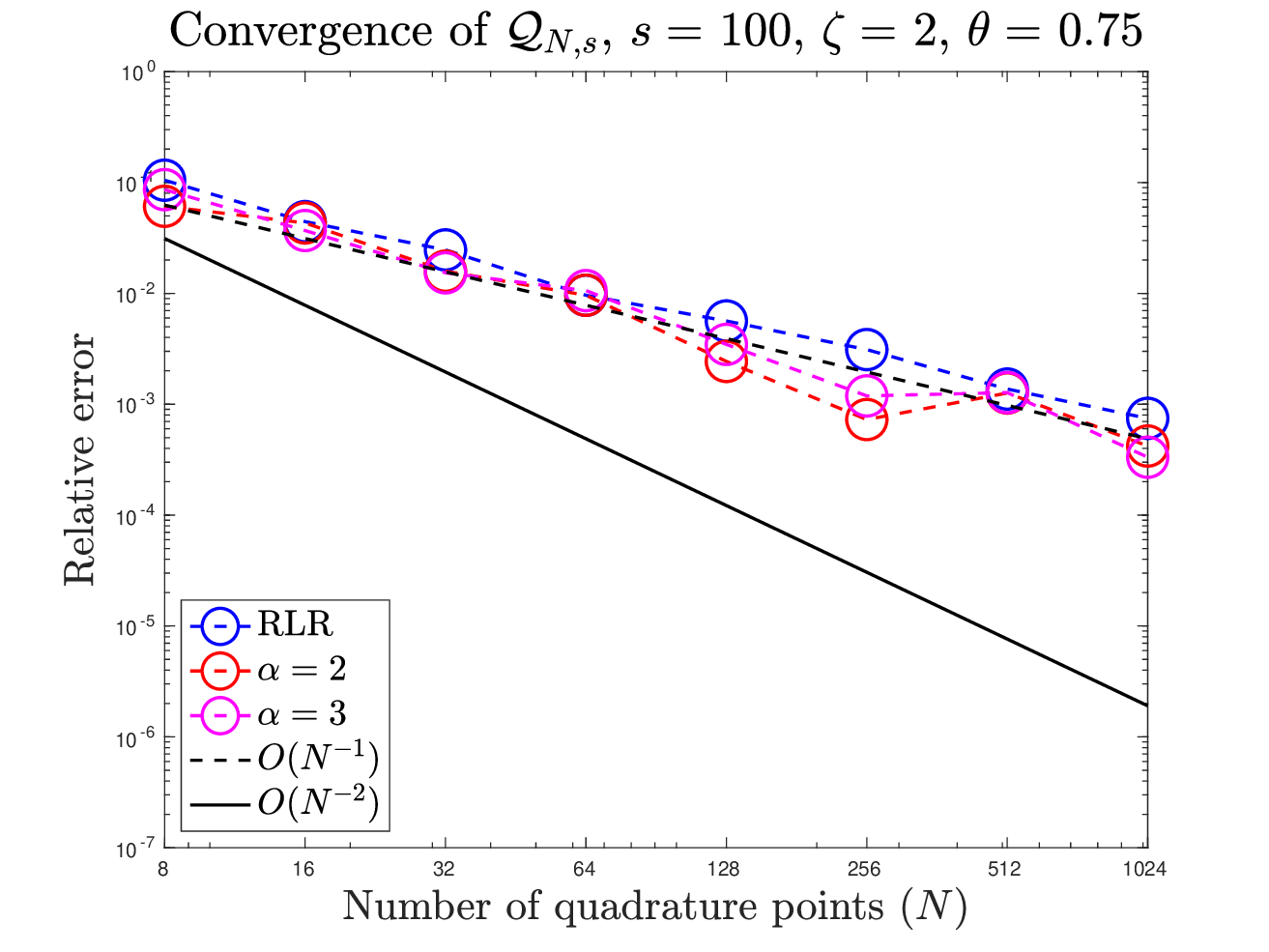}	
		\caption{$\zeta=2, \ s = 100$.}
	\end{subfigure}
	\begin{subfigure}[b]{0.4\textwidth}
		\includegraphics[width=1.0\linewidth]{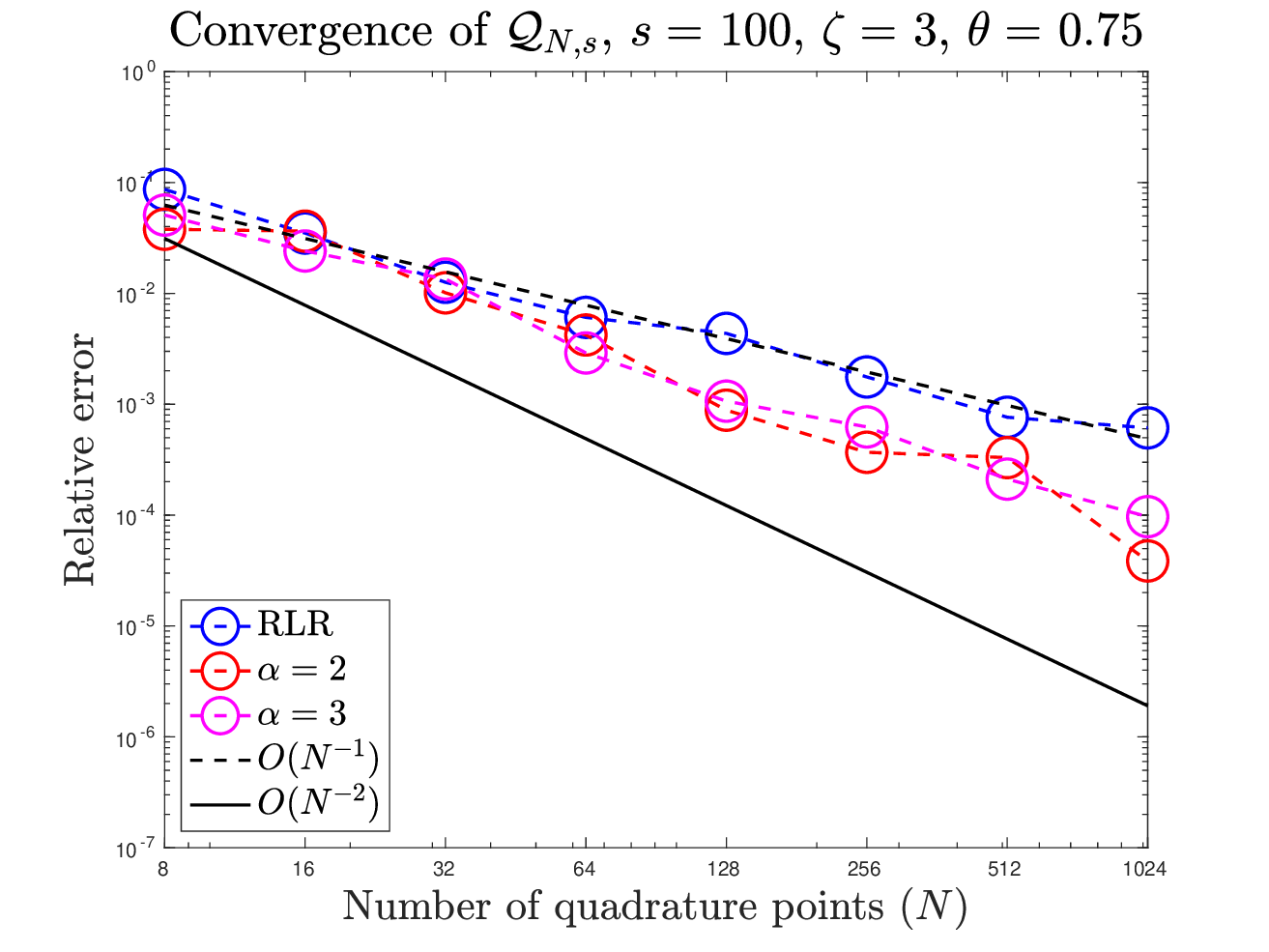}	
		\caption{$\zeta=3, \ s = 100$.}
	\end{subfigure}
	\caption{Convergence of $\mathcal{Q}_{N, s}$ vs. number of QMC points ($N$): $\theta = \tfrac{3}{4}$.}
	\label{fig:Test2-Forward}
\end{figure}

We study the impact of the uncertainty given by the affine-parametric radius \cref{eq:radius} when evaluating the solution at a sample of $K = 10$ equispaced points $\{\bx_k\}_{k=1}^{K}$ in a circle of radius $2$.
We observed the expected order of convergence for the different schemes when $\theta=\tfrac{1}{4}$
(see \Cref{fig:Test1-Forward}). In particular, $N^{1-\zeta}$ is the optimal rate of convergence given the asymptotic behaviour of $\boldsymbol{b}$. For the case of $\theta=\tfrac{3}{4}$ (large deformations) the observed convergence rates are deteriorated as portrayed in \Cref{fig:Test2-Forward}.

\subsection{Inverse UQ: Bayesian Shape Inversion}
We define a ground truth $\y^{\star}\in [-\tfrac{1}{2}, \tfrac{1}{2}]^s$ and compute the solution of \cref{LSvar}. 
Then, we sample measurements at $K = 10$ equispaced points and introduce a noise given by $\eta\in \mathcal{N}(0, \Sigma),$ where $\Sigma$ is a positive definite covariance matrix. In our examples, we simply use
$\Sigma \coloneqq \sigma^2 \mathsf{I}$ with $\sigma = 0.1.$
The measurements $\Upsilon$ are given by \cref{eq:model_noise}
with $\mathcal{G}$ as in \cref{eq:param_to_measurement_map}.

\begin{figure}[ht!]
\centering
	\begin{subfigure}[b]{0.4\textwidth}
		\includegraphics[width=1.0\linewidth]{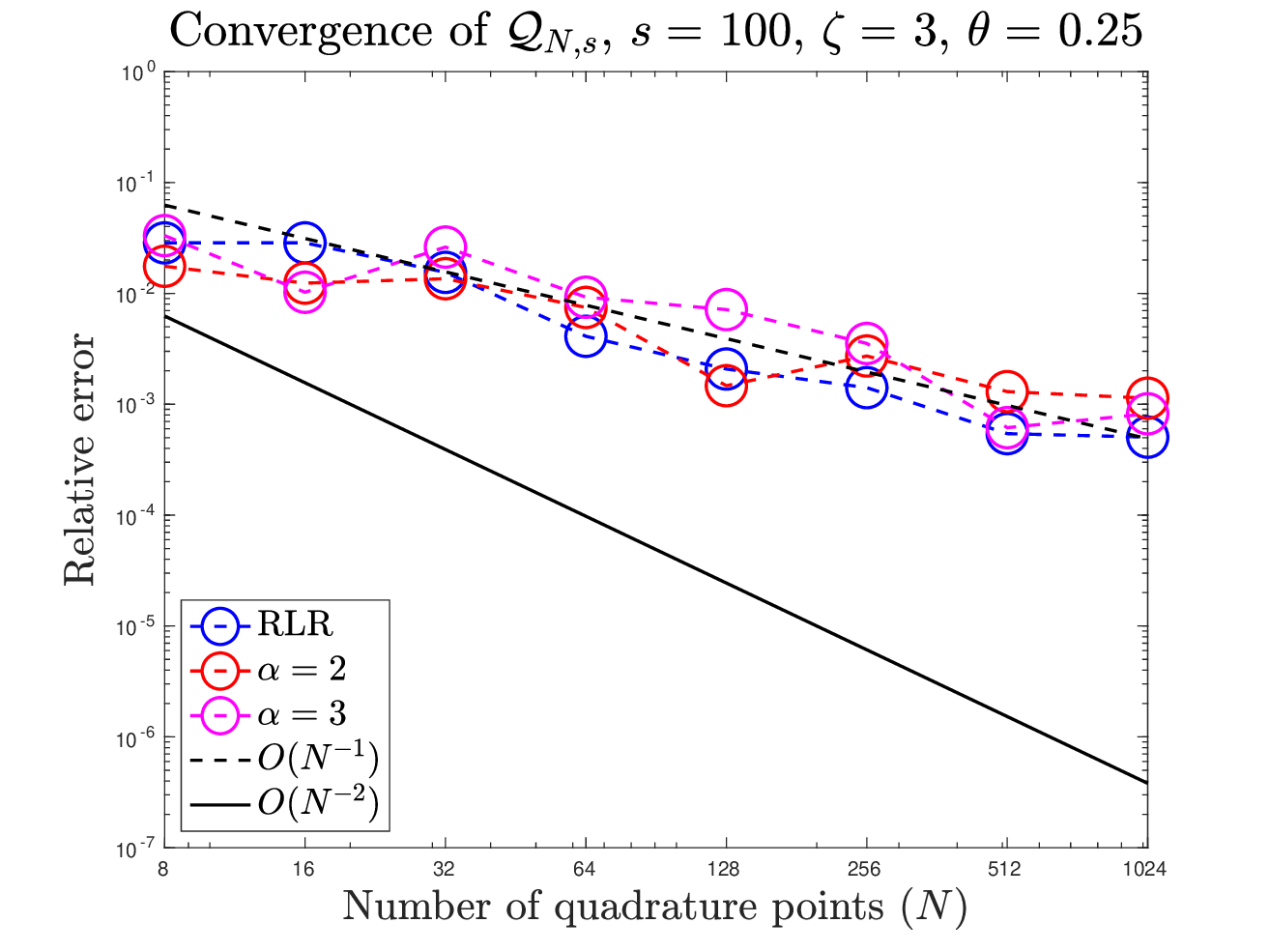}	
		\caption{Convergence of $\mathcal{Q}_{N, s}$.}
        \label{fig:Test1-Inverse-rates}
	\end{subfigure}
	\begin{subfigure}[b]{0.4\textwidth}
         \includegraphics[width=1.0\linewidth]{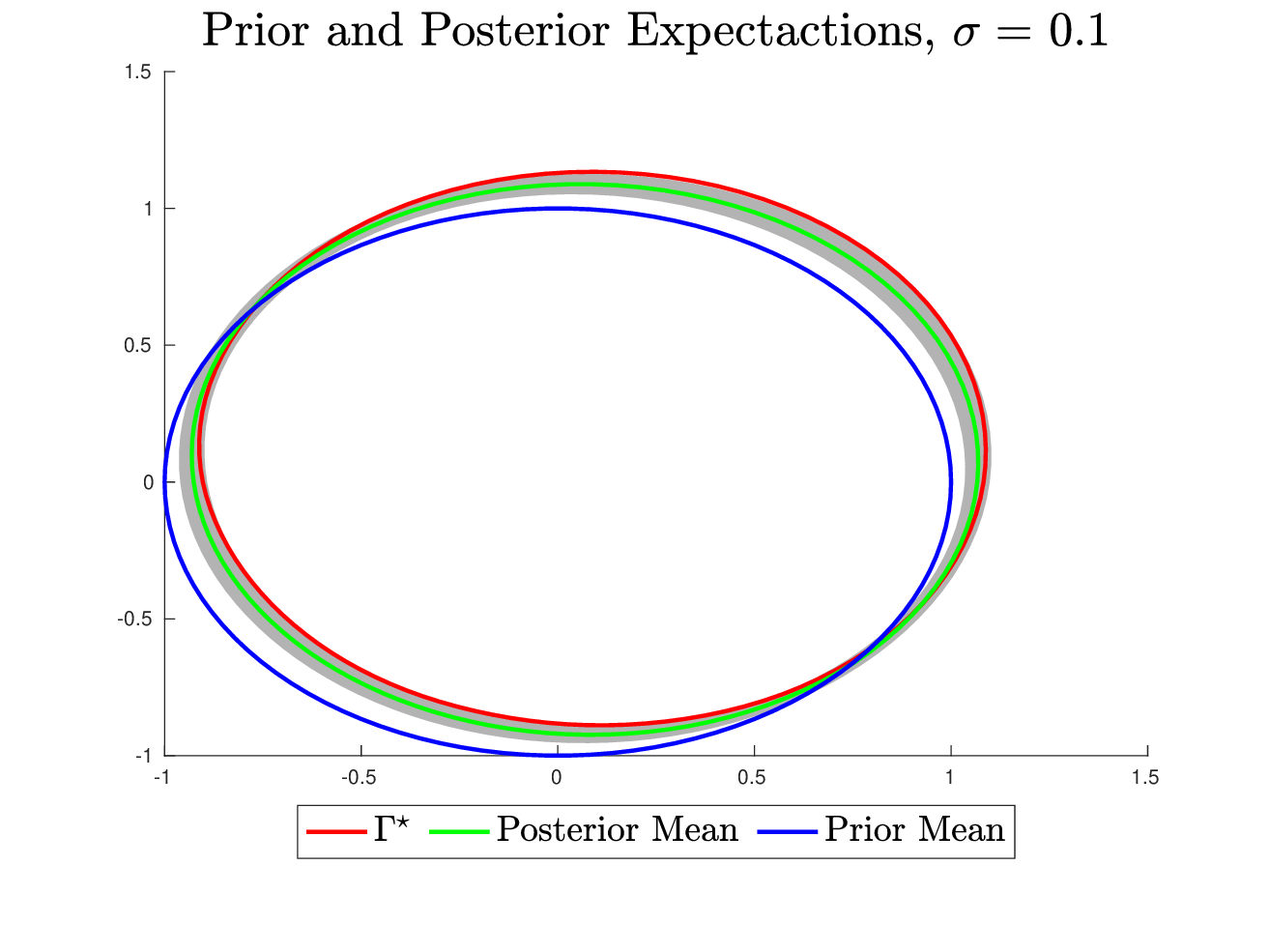}	
		\caption{Bayesian Shape Inversion.}
        \label{fig:Test1-Inverse-disk}
	\end{subfigure}
	\caption{Experiment with $s = 100, \ \zeta=3, \ \theta = \tfrac{1}{4}$ and $\sigma = 0.1$.}
	\label{fig:Test1-Inverse}
\end{figure}

\begin{figure}[ht!]
\centering
	\begin{subfigure}[b]{0.4\textwidth}
        \includegraphics[width=1.0\linewidth]{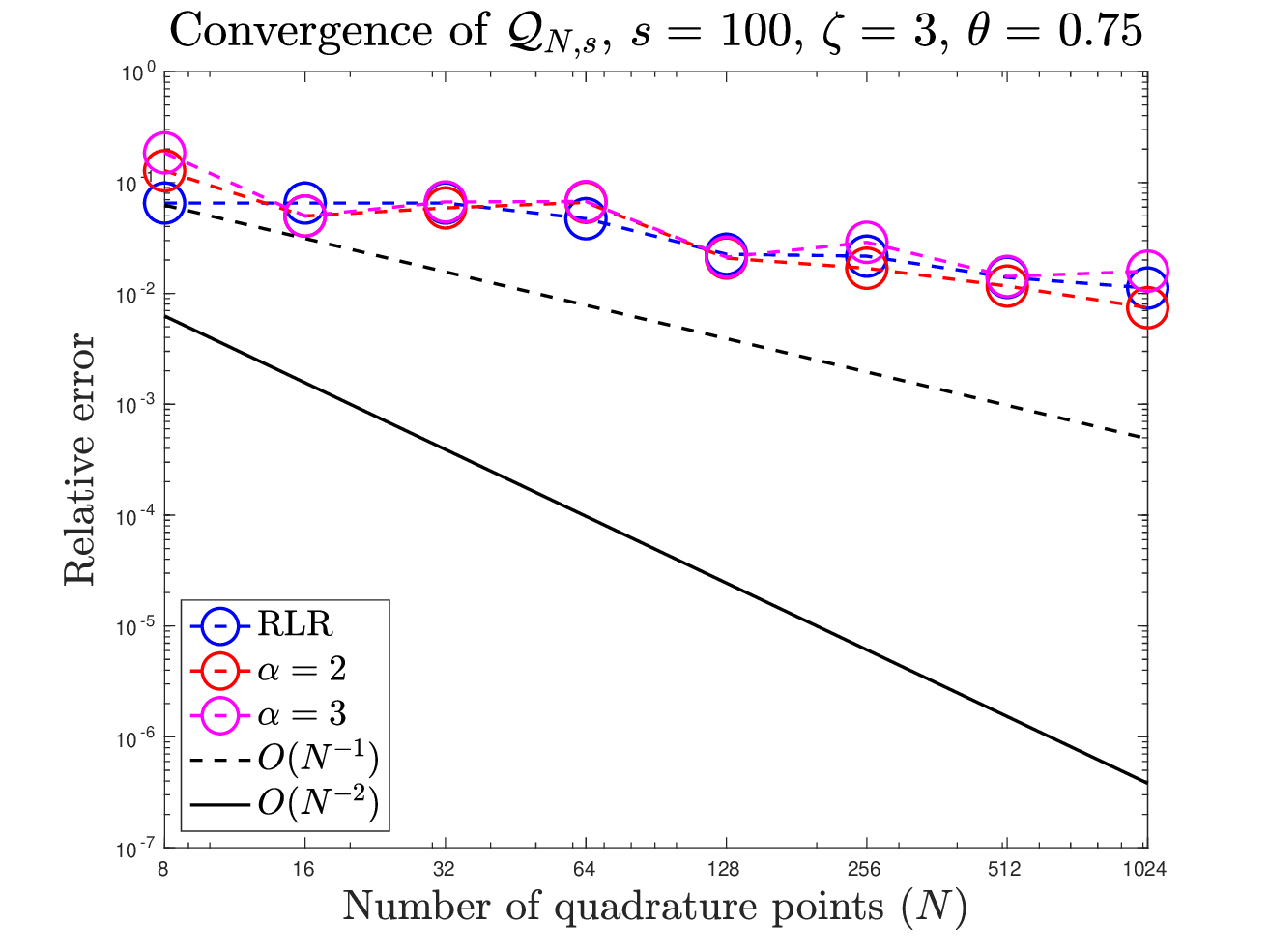}		
		\caption{Convergence of $\mathcal{Q}_{N, s}$.}
        \label{fig:Test2-Inverse-rates}
	\end{subfigure}
	\begin{subfigure}[b]{0.4\textwidth}
		\includegraphics[width=1.0\linewidth]{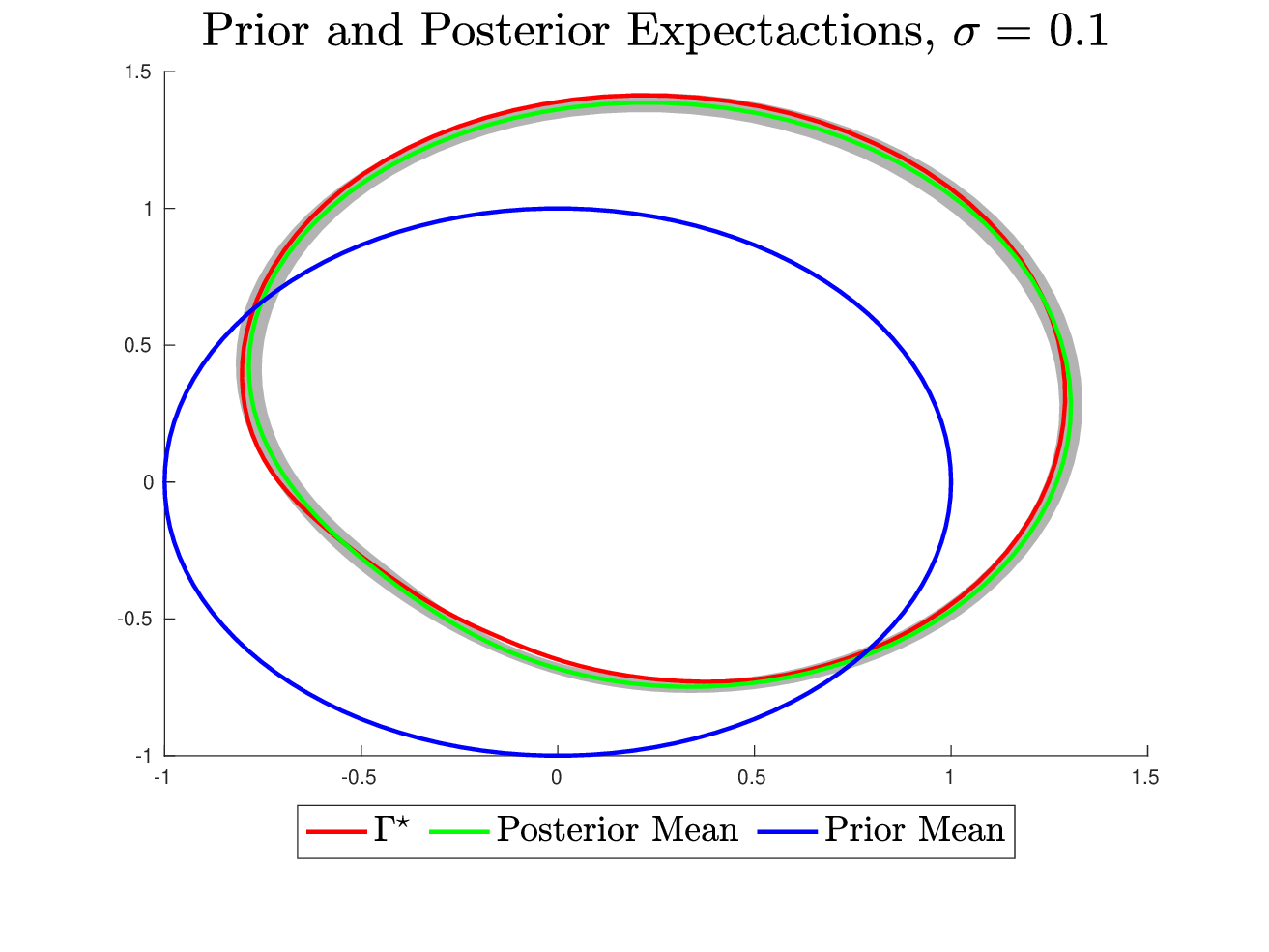}	
		\caption{Bayesian Shape Inversion.}
        \label{fig:Test2-Inverse-disk}
	\end{subfigure}
	\caption{Experiment with $s = 100, \ \zeta=3, \ \theta = \tfrac{3}{4}$ and $\sigma = 0.1$.}
	\label{fig:Test2-Inverse}
\end{figure}

We show convergence results for the expected value of $\phi(\y) = \br_{\y}$ (cf. \Cref{eq:bayesian_shape_iden})
with respect to the posterior distribution (see \cref{fig:Test1-Inverse-rates}), approximated with QMC quadrature rules under the setting described in \Cref{sec:forward-uq}.
We also illustrate the computed prior and posterior expectations compared to the ground truth in Figure \cref{fig:Test1-Inverse-disk}. The same experiment is shown in \cref{fig:Test2-Inverse} for the case of large deformations: $\theta=\tfrac{3}{4}.$ As observed, the difference between prior and posterior expectations becomes relevant in this scenario, despite the limitations given by large deformations in the $(\boldsymbol{b},p,\epsilon)$-holomorphy setting.

\section{Concluding Remarks}
\label{sec:concluding_remarks}
In this work, we propose to tackle the problem of computational UQ for acoustic wave scattering by an inhomogeneous scatterer of uncertain shape. A key component of our approach involves casting this problem, which is posed in an unbounded domain, into a suitable VIE.

Unlike existing methods, such as \cite{HSSS15}, our approach does not require artificial truncation of the computational domain or the imposition of boundary conditions on an artificial boundary to emulate the domain's unboundedness. Traditional approaches to acoustic scattering problems in unbounded domains often use boundary integral formulations. However, this approach is limited when inhomogeneities are not piece-wise constant, as it typically requires a Green's function with an explicit expression for computational efficiency.

It should be noted that using VIOs after a Galerkin-type discretization on the volume results in a dense linear system of equations. This issue also arises with the boundary element method, but in the latter, the Galerkin finite-dimensional subspaces are constructed by meshing the object's boundary rather than the volume. This leads to a larger linear system for the same target accuracy. Recent developments in the Galerkin discretization of VIEs address this problem \cite{LH23,labarca2024coupled}.

A significant milestone in our work is proving that the VIOs depend holomorphically on shape deformations. Consequently, the solution to the well-posed Lippmann-Schwinger equation also depends holomorphically on the domain's perturbations. This property, known as shape holomorphy, is crucial for implementing dimension-robust techniques for computational UQ in problems with high-dimensional, distributed parametric inputs. We present two applications that demonstrate this.

Current and future work includes extending our approach to electromagnetic wave scattering by inhomogeneous inclusions and approximating the parameter-to-solution map using the Galerkin-POD Neural Network method as in \cite{HU18,weder2024galerkin}.



\bibliographystyle{siamplain}
\bibliography{ref-2.bib}

\end{document}